\documentclass[10pt]{article}

\usepackage{amsmath,amsthm,amsfonts,amssymb}
\usepackage{graphicx}

\usepackage[usenames]{color}
\usepackage{hyperref}
\usepackage{soul}

\newenvironment{blue}{\color{blue}}{}

\def\MN{{\mathbb{N}}}
\def\MZ{{\mathbb{Z}}}

\title{Groups with free regular length functions in $\mathbb{Z}^n$}
\author{\textsf{Olga Kharlampovich}
\and \textsf{Alexei Myasnikov}
\and \textsf{Vladimir Remeslennikov}
\and \textsf{Denis Serbin}}
\date{\textsf{}}

\newtheorem{example}{Example}
\newtheorem{corollary}{Corollary}
\newtheorem{prop}{Proposition}
\newtheorem{theorem}{Theorem}
\newtheorem{lemma}{Lemma}
\newtheorem{defn}{Definition}
\newtheorem{remark}{Remark}

\begin{document}
\maketitle

\begin{abstract}

This is the first paper in a series of three where we take on  the unified theory of non-Archimedean group actions, length functions and infinite words.  Our main goal is to show that group actions on $\mathbb{Z}^n$-trees give one a  powerful tool to study groups. All finitely generated  groups acting freely on $\mathbb{R}$-trees also act freely  on some $\mathbb{Z}^n$-trees, but the latter ones form a much larger class.   The natural effectiveness  of all constructions for $\mathbb{Z}^n$-actions (which is not the case for $\mathbb{R}$-trees) comes along with a robust algorithmic theory. In this paper  we describe the algebraic structure of finitely generated groups acting freely and regularly on $\MZ^n$-trees and give  necessary and sufficient conditions for  such actions.

\end{abstract}

\tableofcontents

\section{Introduction}

This is the first paper in a series of three where we take on  the unified theory of non-Archimedean group actions, length functions and infinite words.  Our main goal is to show that group actions on $\mathbb{Z}^n$-trees give one a  powerful tool to study groups. All finitely generated  groups acting freely on $\mathbb{R}$-trees also act freely  on some $\mathbb{Z}^n$-trees, but the latter ones form a much larger class.   The natural effectiveness  of all constructions for $\mathbb{Z}^n$-actions (which is not the case for $\mathbb{R}$-trees) comes along with a robust algorithmic theory. In this paper  we describe the algebraic structure of finitely generated groups acting freely and regularly on $\MZ^n$-trees and give  necessary and sufficient conditions for  such actions.

\medskip \noindent
{\bf Group actions and length functions}.
There are two classical approaches to the fundamental groups of graphs of groups: one uses group actions on trees, and another one - length functions.

In his seminal book \cite{Serre}  Serre laid down fundamentals of the theory of
groups acting freely on simplicial trees. In the following decade
Serre's approach unified several geometric and
combinatorial methods of group theory into a unique tool,
known today as Bass-Serre theory. This was one of the major achievements of combinatorial group theory in 1970's.

In 1980's  Morgan and Shalen  introduced group actions on $\Lambda$-trees for an arbitrary
ordered abelian group $\Lambda$ \cite{MS}. In particular, they studied actions on ${\mathbb R}$-trees in relation with
Thurston's Geometrization Theorem.

Alperin and Bass \cite{AB}  developed the initial framework of the  theory of
 group actions on $\Lambda$-trees and stated the fundamental research goals:
 find the group theoretic information carried by an action (by isometries) on a
$\Lambda$-tree; generalize Bass-Serre theory to actions on arbitrary $\Lambda$-trees.

 A joint effort of several
researchers culminated in a description of finitely generated groups
acting freely on $\mathbb R$-trees \cite{BF,GLP}, which is now known as
Rips' theorem: a finitely generated group acts freely on an
$\mathbb{R}$-tree if and only if it is a free product  of free
abelian groups and surface groups (with exception of non-orientable
groups of genus 1, 2, and 3).
The key ingredient of this theory is the so-called "Rips machine", the idea of which  comes from Makanin's algorithm  for solving equations in free groups \cite{Mak}. The Rips machines appear in applications as a general tool that takes a sequence of isometric actions of a
group $G$ on   some "negatively curved spaces" and produces an isometric
action of $G$ on an $\mathbb R$-tree in the Gromov-Hausdorff limit. Free actions on $\mathbb{R}$-trees  cover all Archimedean actions, since every group
acting freely on a $\Lambda$-tree for an Archimedean ordered abelian group $\Lambda$ acts freely also on an $\mathbb{R}$-tree.

The case of non-Archimedean free  actions is wide open. In \cite{B} Bass studied finitely generated groups acting freely on  ($\Lambda \oplus \mathbb{Z}$)-trees with respect to the lexicographic order on
$\Lambda \oplus \mathbb{Z}$. Recently,  Guirardel (see \cite{G}) studied finitely generated
groups acting freely on an ${\mathbb R}^n$-tree (with the lexicographic order).
 However, the following main problem of the Alperin-Bass program remains largely intact.

\medskip {\bf Problem}
{\it Describe finitely presented (finitely generated) groups acting freely on an arbitrary $\Lambda$-tree.}
\medskip

Much earlier Lyndon introduced groups equipped with  length functions with values in $\Lambda$ as a tool to carry over Nielsen cancelation theory from free groups to a much more general setting \cite{L} (see also \cite{LS}). Some partial results (for $\Lambda = \mathbb{Z}$ or $\Lambda = \mathbb{R}$) were obtained in \cite{Hoare1,Hoare2,Harrison,Prom,AM}. In  \cite{Ch2} Chiswell described a crucial construction that shows that a group with a length function with values in $\mathbb{R}$ comes naturally from some action of the group on an $\mathbb{R}$-tree, and vice versa. It was realized later that a similar  construction holds for arbitrary group with a free Lyndon length function with values in $\Lambda$ \cite{MS}. Thus free group actions and free Lyndon length functions are  just two equivalent languages describing the same objects. We refer to the book \cite{Ch1} for a detailed discussion on the subject.

Nowadays, the geometric method of group actions seems much more in use  than the length functions, even though the Nielsen theory still  gives very powerful results \cite{KW1,KW2,W1,W2,FRS,GFMRS}. There are several reasons for it. Firstly, the Lyndon's abstract axiomatic approach to length functions is less intuitive than group actions. Secondly, the current development of the abstract Nielsen theory is incomplete and insufficient - one of the principal notions of complete actions and regular length functions were not in use until recently (see the discussion below). Thirdly, and this is crucial, an abstract analog of the Rips machine for length functions with values in $\mathbb{R}$ (or an arbitrary $\Lambda$) was not introduced or  developed enough to be robust in applications. Notice, that in case of $\Lambda = \mathbb{Z}$ the completeness and regularity come almost for free, so there was no need  for axiomatic   formalization, -  it went mostly unnoticed. The regularity axiom  appeared first in \cite{MR,MRS} as a tool to deal with length functions in $\mathbb{Z}^n$ (with respect to the lexicographic order). On the other hand, an analog of the Rips machine for $\Lambda = \mathbb{Z}$ did  exist for a long time  -  the original  Makanin-Razborov process \cite{Raz} for solving equations in free groups. But it was not recognized as such  until \cite{KM1,KM2}, where it was used systematically to get splittings of groups.

\medskip\noindent
{\bf Motivation and recent progress}.
Introduction of infinite $\Lambda$-words  was one of the major recent developments   in the theory of group actions (see Section \ref{subs:2.1} for definitions). In  \cite{MRS}   Myasnikov,
Remeslennikov  and Serbin showed that groups admitting faithful representations by $\Lambda$-words act freely on some $\Lambda$-trees, while Chiswell proved the converse \cite{Ch}. This gives another equivalent approach to group actions. Now  one can bypass the axiomatic view-point on length functions and work instead  with $\Lambda$-words in the same manner as with the ordinary words in the standard free groups, which is much more intuitive and practical. In particular,  one can deal with subgroups of groups acting freely on $\mathbb{Z}^n$-trees (for example, fully residually free groups) using the standard Stallings' folding  argument \cite{MRS,MRS2,KMRS,NS}.

Formulation of the regularity axiom, discussed above,  is another important ingredient of the theory of group actions. This axiom ensures that the group is complete with respect to the Gromov's inner product, i.e., the Gromov's  product of two elements of the  group  can be realized as the length of a particular element of the group (not only as an abstract length).  In the language of actions this means
that  all branch points of the group action are in
one same orbit, as well as the base point. This allows one to use the Nielsen cancelation argument at full strength, very similar to the case of a  free group.  The regularity (completeness) condition is crucial  for  existence of Makanin-Razborov's type   processes  over groups with $\Lambda$-length functions.  We discuss this in the third paper of the series for $\Lambda= \mathbb{Z}^n$, while in the second one we show that every group acting freely on a $\mathbb{Z}^n$-tree isometrically embeds into a group with a regular action on some $\mathbb{Z}^m$-tree.

\medskip\noindent
{\bf Results}. In Section \ref{sec:prelim} we recall the notions of length functions, infinite words, and discuss the regularity axiom.  Here we provide  some crucial examples of groups with regular length functions. In particular, we show that every finitely generated group acting freely on an $\mathbb{R}$-tree  has a  free regular length function in $\mathbb{Z}^n$. More generally, we show that certain  types of HNN extensions have free regular length functions that extend a given  regular free length function on the base groups. Notice, that these are precisely the HNN extension  that play the central role in the elimination process from \cite{KM1}.

Commutation of $\mathbb{Z}^n$-words is an important technical instrument, which can be viewed as a "non-standard" version of the commutation in free groups, we discuss it  in Section \ref{sec:commutation}.

In Section \ref{sec:main} we   describe the algebraic structure  of a finitely generated group $G$ with a regular length function in $\mathbb{Z}^n$ in
terms of HNN-extensions of a very particular type. To this end in
Section \ref{subsec:elem-moves} we introduce Nielsen-like moves and
show that every finite set of generators of $G$ can be transformed
by a sequence of elementary moves into  a minimal one (relative to the $\Lambda$-length). Then, in Section \ref{subsec:minimal}, we show that any
minimal finite set of generators  of $G$ gives an HNN-splitting of $G$.
Here we use the technique of Stallings' pregroups and their universal groups \cite{S}. This is an analog of the celebrated Nielsen theorem on subgroups of free groups. Interestingly, the whole approach closely resembles the Gr${\rm \ddot{o}}$bner basis  method and Buchberger's  algorithm.  Theorem \ref{th:main2}  in Section \ref{subsec:structure} gives the description of the algebraic structure of $G$ in terms of a finite sequence of HNN extensions of a certain type.
In Section \ref{sec:embedding} we prove that the converse of Theorem \ref{th:main2}  also holds, thus finishing the complete description of groups acting freely and regularly on $\mathbb{Z}^n$-trees.

\section{Length functions, actions, infinite words}
\label{sec:prelim}

Here we introduce basic definitions and notations which are to be
used throughout the whole paper.

\subsection{Lyndon length functions and free actions}
\label{subsec:Lyndon}

Let $G$ be a group and $A$ an ordered
abelian group. Then a function $l: G \rightarrow A$ is called a {\it
(Lyndon) length function} on $G$ if the following conditions hold:
\begin{enumerate}
\item [(L1)] $\forall\ g \in G:\ l(g) \geqslant 0$ and $l(1)=0$;
\item [(L2)] $\forall\ g \in G:\ l(g) = l(g^{-1})$;
\item [(L3)] $\forall\ g, f, h \in G:\ c(g,f) > c(g,h)
\rightarrow c(g,h) = c(f,h)$,

\noindent where $c(g,f) = \frac{1}{2}(l(g)+l(f)-l(g^{-1}f))$.
\end{enumerate}
Notice that  $c(g,f)$ may not be defined in $A$ (if
$l(g)+l(f)-l(g^{-1}f)$ is not divisible by $2$), so in the axiom
(L3) we assume that  $A$ is canonically embedded into a divisible
ordered abelian group  $A_{\mathbb{Q}} = A \otimes_\mathbb{Z} \mathbb{Q}$ (see \cite{MRS} for
details).

It is not difficult to derive the following two properties of
length functions from the axioms (L1)-(L3):
\begin{itemize}
\item $\forall\ g, f \in G:\ l(g f) \leqslant l(g) + l(f)$; \item
$\forall\ g, f \in G:\ 0 \leqslant c(g,f) \leqslant
min\{l(g),l(f)\}$.
\end{itemize}
A length function $l:G \rightarrow A$ is called {\em free} if it satisfies the following two axioms.

\begin{enumerate}
\item [(L4)] $\forall\ g, f \in G:\ c(g,f) \in A.$
\item [(L5)] $\forall\ g \in G:\ g \neq 1 \rightarrow l(g^2) > l(g).$
\end{enumerate}

\begin{theorem} \cite{Ch1}
A group $G$ has a free Lyndon length function with values in $A$ if
and only if $G$ acts freely on an $A$-tree.
\end{theorem}
\smallskip

For elements $g_1, \ldots, g_n \in G$ we write
$$g = g_1 \circ \cdots \circ g_n$$
if $g = g_1 \cdots g_n$ and $l(g) = l(g_1) + \cdots + l(g_n)$.
Also, for $\alpha \in A$ we write $g = g_1 \circ_\alpha g_2$ if $g
= g_1 g_2$ and $c(g_1^{-1},g_2) < \alpha$.

\subsection{Regular length functions}
\label{subsec:Lyndon-actions}

In this section we define regular length functions,
and show some examples of groups with regular length functions.
Theorem \ref{main4} gives a host of new examples of such groups.
In the subsequent paper \cite{KMS2} we discuss a geometric
characterization of groups with regular length functions.

A length function $l: G \rightarrow A$ is  called {\it regular} if
it satisfies the {\it regularity} axiom:
\begin{enumerate}
\item [(L6)] $\forall\ g, f \in G,\ \exists\ u, g_1, f_1 \in G:$
$$g = u \circ g_1 \ \& \  f = u \circ f_1 \ \& \ l(u) = c(g,f).$$
\end{enumerate}

Here are several examples of groups with regular free length functions.

\begin{example}
\label{ex:1}
Let $F = F(X)$ be a free group on $X$. The length function
$$|\ | : F \rightarrow \mathbb{Z},$$
where $|f|$ is a the length  of $f \in F$ as a word in $X^{\pm 1}$, is regular.
\end{example}

The following is a more general example.

\begin{example} \label{ex:1b}
Let $F = F(X)$ be a free group on $X$, $H$ a finitely generated subgroup of $H$, and $l_H$ the restriction to $H$ of the length function in $F$ relative to $X$. Then $l_H$ is a regular length function on $H$ if and only if there exists a basis $U$ of $H$ such that every two non-equal elements from $U^{\pm 1}$ have different initial letters.
\end{example}
\begin{proof}
The straightforward proof is not very hard. On the other hand, one can argue using the transformations $\mu$ and $\nu$ from Section  \ref{subsec:elem-moves}. Indeed, any reduced generating set $U$ of $H$ satisfies the conditions above, since the transformations $\mu$ and $\nu$ do not reduce the total length. Conversely, if a generating set $U$ satisfies the conditions then  there is no cancelation in the products $uv$ for $u, v \in U^{\pm 1}$, $u \neq v$, and the result follows as in Example \ref{ex:1}.
\end{proof}

\begin{example}
\label{ex:2}  \cite{MRS}
Lyndon's free $\mathbb{Z}[t]$-group $F^{\mathbb{Z}[t]}$ has a regular
free length function with values in $\mathbb{Z}[t]$.
\end{example}

The next example is just a particular case of Theorem \ref{main4} below.

\begin{example}
\label{ex:3}
Let $F = F(X)$ be a free group with basis $X$, $|\ |$ the standard length function on $F$ relative to $X$, and  $u,v \in F$  such that $|u| = |v|$ and $u$ is not conjugate to $v^{-1}$. Then the  HNN-extension
$$G = \langle F, s \mid u^s = v \rangle,$$
has  a regular free length function $l : G \rightarrow \mathbb{Z}^2$
which extends $|\cdot|$.
\end{example}

\begin{example}
\label{ex:4} For any $n \geq 1$ the orientable surface   group
$$G = \langle x_1, x_2, \ldots, x_{2n-1}, x_{2n} \mid [x_1,x_2] \cdots
[x_{2n-1}, x_{2n}] = 1 \rangle$$
has  a regular free length function $l : G \rightarrow \mathbb{Z}^2$.
\end{example}

\begin{proof}
It suffices to represent  $G$ as an HNN extension from Example \ref{ex:3}.  The word
$$R(X) = x_1 \cdots x_{2n} x_1^{-1} \cdots x_{2n}^{-1}.$$
is quadratic, so   there exists an automorphism $\phi$ of $F =
F(x_1,\ldots, x_{2n})$ such that
$$R(X)^\phi = [x_1,x_2] \cdots [x_{2n-1}, x_{2n}]$$
(see, for example,   Proposition 7.6
\cite{LS}).  It follows that $G$ is isomorphic to
$$G^\prime = \langle x_1,\ldots, x_{2n} \mid x_1 \cdots x_{2n} x_1^{-1}
\cdots x_{2n}^{-1} = 1 \rangle,$$
which can be represented as an HNN-extension of the required form
$$G^\prime = \langle F(x_2,\ldots, x_{2n}), x_1 \mid x_1 (x_2 \cdots x_{2n})
x_1^{-1} = x_{2n} x_{2n-1} \cdots x_2 \rangle,$$ since $|x_2
\cdots x_{2n}| = |x_{2n} x_{2n-1} \cdots x_2|$.

\end{proof}

\begin{example}
\label{ex:5} For any $n,\
n\geq 3$ the non-orientable surface group
$$G = \langle x_1, x_2, \ldots,  x_{n} \mid x_1^2x_2^2\ldots x_n^2 = 1\rangle $$
has a regular free length function $l : G \rightarrow \mathbb{Z}^2$.
 \end{example}
 \begin{proof}
 Again, it suffices to represent $G$ as an HNN extension from Example \ref{ex:3}.

 An argument similar to the one in the proof of Example \ref{ex:4} shows that the group $G$ is isomorphic to
  $$G^\prime = \langle x_1, x_2, \ldots,  x_{n} \mid x_1\ldots x_{n-1}x_nx_1^{-1}\ldots
x_{n-1}^{-1}x_n \rangle $$
  and the result follows, since the presentation above can be written as
 $$G^\prime = \langle x_1, x_2, \ldots,  x_{n} \mid x_1 (x_2\ldots x_{n-1}x_n)x_1^{-1} = x_n^{-1}x_{n-1}\ldots x_2 \rangle $$

\end{proof}

\begin{example}
\label{ex:6} A free abelian group of rank $n$ has a free regular
length function in ${\mathbb Z}^n$.
\end{example}
\begin{proof}
Let $G = \mathbb{Z}^n$ be a free abelian group of rank $n$.   Then $G$ is an ordered abelian group relative to the right lexicographic order $\leq$. The absolute value $|u|$ of an element $u \in G$, defined as $|u| = \max \{u,-u\}$,  gives a
free length function $l : G \rightarrow \mathbb{Z}^n$.  It is easy to see that $l$ is regular.
\end{proof}

\begin{example}
\label{ex:7} Let  $G_i$, $i = 1,2$, be a group  with a free
regular length function $l_i : G_i \to {\mathbb Z}^n$. Then the free product $G = G_1\ast G_2$  has a free regular
length function in $l:G \to {\mathbb Z}^n$ that extends the functions $l_1$ and $l_2$.
\end{example}
 \begin{proof}
Let $g \in G$ given in the reduced form  $g=u_1v_1\ldots u_kv_k,$ where $u_1,\ldots, u_k\in G_1$ and
$v_1,\ldots, v_k\in G_2$. Define $l:G \to  {\mathbb Z}^n$ by
$$l(g)=\sum _{i=1}^k (l_1(u_i)+l_2(v_i).$$
It is not hard to see that $l$ is a regular free length function (see, for example, \cite{LS}).

\end{proof}

\begin{example} Let  $G$ be a finitely generated  group  acting freely on $\mathbb{R}$-tree.
 Then $G$ has a free regular length function in ${\mathbb
Z}^n$, where $n$ is the maximal rank of free abelian
subgroups (centralizers) of $G$.
\end{example}
 \begin{proof}
 By Rips, theorem every finitely generated group acting freely on an $\mathbb{R}$-tree is a free product of groups described in Examples \ref{ex:4}, \ref{ex:5}, \ref{ex:6}, hence the result.
  \end{proof}

\subsection{Infinite words and length functions}
\label{subs:2.1}

In this subsection at first we recall some notions from the theory
of ordered abelian groups (for all the details refer to the books
\cite{Gl} and \cite{KopMed}) and then following \cite{MRS}
describe the construction of infinite words.

\smallskip

Let $A$ be an ordered abelian group. Recall that for elements $a,b
\in A$ the {\it closed segment} $[a,b]$ is defined by
$$[a,b] = \{c \in A \mid a \leqslant c \leqslant b \}.$$
Now a  subset $C \subset A$ is called {\it convex} if for every $a,
b \in C$ the set $C$ contains  $[a,b]$. In particular, a subgroup
$B$ of $A$ is convex if $[0,b] \subset B$ for every positive $b
\in B$.

An ordered abelian group $A$ is called {\it discretely ordered} if
$A$ has a minimal positive element (we denote it by $1_A$). In
this event, for any $a \in A$ the following hold:
\begin{enumerate}
\item[1)] $a+1_A = \min\{b \mid b > a\}$,
\item[2)] $a-1_A = \max\{b \mid b < a\}$.
\end{enumerate}

Observe that if $A$ is any ordered abelian group then $\mathbb{Z}
\oplus A$ is discretely ordered with respect to the right
lexicographic order.

\smallskip

Let $A$ be a discretely ordered abelian group and let $X = \{x_i
\mid i \in I\}$ be a set. Put $X^{-1} = \{x_i^{-1} \mid i \in I\}$
and $X^\pm = X \cup X^{-1}$. An {\em $A$-word} is a function of
the type
$$w: [1_A,\alpha_w] \to X^\pm,$$
where $\alpha_w \in A,\ \alpha_w \geqslant 0$. The element
$\alpha_w$ is called the {\em length} $|w|$ of $w$.

\smallskip

By $W(A,X)$ we denote the set of all $A$-words. Observe, that
$W(A,X)$ contains an empty $A$-word which we denote by
$\varepsilon$.

Concatenation $uv$ of two $A$-words $u,v \in W(A,X)$ is an
$A$-word of length $|u| + |v|$ and such that:
\[ (uv)(a) = \left\{ \begin{array}{ll}
 \mbox{$u(a)$}  & \mbox{if $1_A \leqslant a \leqslant |u|$} \\
 \mbox{$v(a - |u|)$ } & \mbox{if $|u|  < a \leqslant |u| + |v|$}
\end{array}
\right. \]

An $A$-word $w$ is {\it reduced} if $w(\beta + 1_A) \neq
w(\beta)^{-1}$ for each $1_A \leqslant \beta < |w|$. We denote by
$R(A,X)$ the set of all reduced $A$-words. Clearly, $\varepsilon
\in R(A,X)$.

\smallskip

For $u \in W(A,X)$ and $\beta \in [1_A, \alpha_u]$ by $u_\beta$ we
denote the restriction of $u$ on $[1_A,\beta]$. If $u \in R(A,X)$
and $\beta \in [1_A, \alpha_u]$ then
$$u = u_\beta \circ {\tilde u}_\beta,$$
for some uniquely defined ${\tilde u}_\beta$.

An element ${\rm com}(u,v) \in R(A,X)$ is called the
({\emph{longest}) {\it common initial segment} of $A$-words $u$
and $v$ if
$$u = {\rm com}(u,v) \circ \tilde{u}, \ \ v = {\rm com}(u,v) \circ
\tilde{v}$$ for some (uniquely defined) $A$-words $\tilde{u},
\tilde{v}$ such that $\tilde{u}(1_A) \neq \tilde{v}(1_A)$.

Now, we can define the product of two $A$-words. Let $u,v \in
R(A,X)$. If ${\rm com}(u^{-1}, v)$ is defined then
$$u^{-1} = {\rm com}(u^{-1},v) \circ {\tilde u}, \ \ v = {\rm com}
(u^{-1},v) \circ {\tilde v},$$ for some uniquely defined ${\tilde
u}$ and ${\tilde v}$. In this event put
$$u \ast v = {\tilde u}^{-1} \circ {\tilde v}.$$
The  product ${\ast}$ is a partial binary operation on $R(A,X)$.

\smallskip

An element $v \in R(A,X)$ is termed {\it cyclically reduced} if
$v(1_A)^{-1} \neq v(|v|)$. We say that an element $v \in R(A,X)$
admits a {\it cyclic decomposition} if $v = c^{-1} \circ u \circ
c$, where $c, u \in R(A,X)$ and $u$ is cyclically reduced. Observe
that a cyclic decomposition is unique (whenever it exists). We
denote by $CR(A,X)$ the set of all cyclically reduced words in
$R(A,X)$ and by $CDR(A,X)$ the set of all words from $R(A,X)$
which admit a cyclic decomposition.

\smallskip

Below we refer to $A$-words as {\it infinite words} usually
omitting $A$ whenever it does not produce any ambiguity.

The following result establishes the connection between infinite
words and length functions.
\begin{theorem}
\label{co:3.1} \cite{MRS} Let $A$ be a discretely ordered abelian
group and $X$ be a set. Then any subgroup $G$ of $CDR(A,X)$ has a
free Lyndon length function with values in $A$ -- the restriction
$L|_G$ on $G$ of the standard length function $L$ on $CDR(A,X)$.
\end{theorem}

The converse of the theorem above was obtained by I.Chiswell
\cite{Ch}.

\begin{theorem}
\label{chis} \cite{Ch} Let $G$ have a free Lyndon length function
$L : G \rightarrow A$, where $A$ is a discretely ordered abelian
group. Then there exists a set $X$ and a length preserving embedding
$\phi : G \rightarrow CDR(A,X)$, that is, $|\phi(g)| = L(g)$ for any
$g \in G$.
\end{theorem}

\begin{corollary}
\label{chis-cor} \cite{Ch} Let $G$ have a free Lyndon length
function $L : G \rightarrow A$, where $A$ is an arbitrary ordered
abelian group. Then there exists an embedding $\phi : G
\rightarrow CDR(A',X)$, where $A' = \mathbb{Z} \oplus A$ is
discretely ordered with respect to the right lexicographic order
and $X$ is some set, such that, $|\phi(g)| = (0,L(g))$ for any $g
\in G$.
\end{corollary}

Theorems \ref{co:3.1} and \ref{chis}, and Corollary \ref{chis-cor}
show that a group has a free Lyndon length function if and only if
it embeds into a set of infinite words and this embedding
preserves the length. Moreover, it is not hard to show that this
embedding also preserves regularity of the length function.

\begin{theorem}
\label{chis-cor-1}
\cite{KMS} Let $G$ have a free regular Lyndon length
function $L : G \rightarrow A$, where $A$ is an arbitrary ordered
abelian group. Then there exists an embedding $\phi : G
\rightarrow R(A',X)$, where $A'$ is a discretely ordered abelian
group and $X$ is some set, such that, the Lyndon length function
on $\phi(G)$ induced from $R(A',X)$ is regular.
\end{theorem}

\section{Commutation in infinite words}
\label{sec:commutation}

Let $G$ be a subgroup of $CDR(A,X)$, where $A$ is a discretely ordered
abelian group. We fix $G$ for the rest of the section.

It is not hard to see that the set $\{A_i \mid i \in I_A\}$ of
convex subgroups of $A$ is linearly
ordered by inclusion, that is, for any $i,j \in I_A,\ i \neq j$ we
have $A_i < A_j$ whenever $i < j$ and
$$A = \bigcup_{i\in I_A} A_i.$$
We say that $g \in G$ {\em has the height} $i \in I_A$ and denote
$ht(g) = i$ if $|g| \in A_i$ and $|g| \notin A_j$ for any $j < i$.
Observe that this definition depends only on $G$ since the complete
chain of convex subgroups of $A$ is unique.

It is easy to see that
$$ht(g_1 g_2) \leq \max\{ht(g_1),ht(g_2)\},$$
hence, if $G = \langle g_1, \ldots, g_k \rangle$ then we define
$$ht(G) = \max\{ht(g_1),\ldots,ht(g_k)\}.$$

Using the characteristics of elements of $G$ introduced above
we prove several technical results.

\begin{lemma}
\label{le:LS} Let $f,h \in G$ be cyclically reduced. If\/ $c(f^m, h^n)
\geqslant |f| + |h|$ for some $m, n > 0$ then $[f,h] = \varepsilon$.
\end{lemma}
\begin{proof} Suppose $|h| \geqslant |f|$ and $c(f^m, h^n) \geqslant
|f| + |h|$ for some $m, n > 0$. Notice that $c(f^m, h^n) \geqslant |f| + |h|$
implies $ht(f) = ht(h)$. Hence, there exists $k \in \mathbb{N}$ such that
$|h| \geq k|f|,\ |h| \leq (k+1)|f|$.

We have $h = f^k \circ h_1,\ |f| > |h_1|,\ h_1 \in G,\ k \geqslant 1$ and $f =
h_1 \circ f_1$. Since $c(f^m, h^n) \geqslant |f| + |h|$ one has
$(f^k \circ h_1) \circ f = f^{k+1} \circ h_1$. So, $h_1 \circ h_1
\circ f_1 = h_1 \circ f_1 \circ h_1$ and $f = h_1 \circ f_1 = f_1
\circ h_1$. It follows that $[f_1,h_1] = \varepsilon$, hence,
$[h_1, f] = \varepsilon$, and $[f,h] = \varepsilon$.
\end{proof}

\begin{lemma}
\label{le:cycl}
For any two $g_1,g_2 \in G$ if $[g_1,g_2] = \varepsilon$ and
$g_1 = c^{-1} \circ h_1 \circ c,\ g_2 = d^{-1} \circ h_2 \circ d$
are their cyclic decompositions then $c = d$.
\end{lemma}
\begin{proof} Without loss of generality we can assume that $c = \varepsilon$.
From $[g_1,g_2] = \varepsilon$ we get
$$(d^{-1} \circ h_2^{-1} \circ d) \ast h_1 \ast (d^{-1} \circ h_2
\circ d) = h_1,$$
where $h_1 \neq \varepsilon,\ h_2 \neq \varepsilon$. Since $h_1$ is cyclically
reduced, either
$$(d^{-1} \circ h_2^{-1} \circ d) \ast h_1 = (d^{-1} \circ h_2^{-1}
\circ d) \circ h_1,$$
or
$$h_1 \ast (d^{-1} \circ h_2 \circ d) = h_1 \circ (d^{-1} \circ h_2
\circ d).$$
Assume the latter (otherwise the argument is similar). Hence, we have
$$h_1 \circ (d^{-1} \circ h_2^{-1} \circ d) = (d^{-1} \circ h_2^{-1}
\circ d) \circ  h_1,$$
and unless $d = \varepsilon$ we have $c(d^{-1},h_1^{-1}) > 0$, that is,
$h_1 \ast d^{-1} \neq h_1 \circ d^{-1}$ - a contradiction to our assumption.
Thus, $d = \varepsilon$.

\end{proof}

In particular, it follows that if $g \in G$ is cyclically reduced then all
elements of $C_G(g)$ are cyclically reduced as well. In this event we call $C_G(g)$
{\it cyclically reduced}.

\begin{lemma}
\label{le:0} Let $f,h \in G$ be such that $h$ is cyclically reduced and $ht(f) > ht(h)$.
If $ht(f^{-1} \ast h \ast f) < ht(f)$ then for every $n \in \mathbb{N}$ either
$f = h^n \circ f_n$, or $f = h^{-n} \circ f_n$ for some $f_n \in G$.
\end{lemma}
\begin{proof} Since $h$ is cyclically reduced, either $f^{-1} \ast h = f^{-1} \circ h$ or
$h \ast f = h \circ f$. Assume the latter (otherwise the argument is similar).

Assume that the statement holds for some $k \geq 0$, that is, $f = h^k \circ f_k$.
At the same time, $f^{-1} \ast (h \circ f) = f_k^{-1} \ast (h \circ f_k)$ and
$ht(f_k^{-1} \ast (h \circ f_k)) < ht(f_k) = ht(f)$ implies that $h$ cancels completely in
the product $f_k^{-1} \ast (h \circ f_k)$. Hence, $f_k$ has $h$ as an
initial subword and the statement holds for $k + 1$.

\end{proof}

\begin{lemma}
\label{le:1} Let $f,h_1,h_2 \in G$ be such that $ht(h_1),ht(h_2) < ht(f)$ and
$ht(f^{-1} \ast h_1 \ast f),\ ht(f^{-1} \ast h_2 \ast f) < ht(f)$.
Then $[h_1,h_2] = \varepsilon$.
\end{lemma}
\begin{proof} Without loss of generality we can assume $h_1$ to be cyclically reduced.
Indeed,
if $h_1 = d^{-1} \circ \overline{h_1} \circ d$ is the cyclic decomposition of $h_1$
then it follows that $f = d^{-1} \circ g$ and
$$f^{-1} \ast h_1 \ast f = g^{-1} \ast \overline{h_1} \ast g,\ f^{-1}
\ast h_2 \ast f = g^{-1} \ast (d \ast \overline{h_2} \ast d^{-1}) \ast g,$$
so, one can consider the triple $g,\overline{h_1},d \ast \overline{h_2} \ast d^{-1}$
instead of $f,h_1,h_2$.

Since $h_1$ is cyclically reduced it follows that in the product $f^{-1} \ast h_1 \ast f$
either $f^{-1} \ast h_1 = f^{-1} \circ h_1$ or $h_1 \ast f = h_1 \circ f$. Assume the latter.
Now, from $ht(f^{-1} \ast h_1 \ast f) < ht(f)$ by Lemma \ref{le:0}, $f$ contains any
positive power of $h_1$ as an initial subword.

At the same time, $h_2 = c^{-1} \circ \overline{h_2} \circ c$, where $\overline{h_2}$ is
cyclically reduced, and $f = c^{-1} \circ f_1$. Thus, $f^{-1} \ast h_2 \ast f =
f_1^{-1} \ast \overline{h_2} \ast f_1,\ ht(f_1^{-1} \ast \overline{h_2} \ast f_1) < ht(f_1)$
and by Lemma \ref{le:0} we can assume that $f_1$ contains any positive power of
$\overline{h_2}$ as an initial segment (for negative power the argument is similar).

Consider several cases.

\begin{enumerate}

\item [{\bf (1)}] $ht(c) \leqslant ht(h_1)$

\begin{enumerate}

\item[{\bf (a)}] $ht(h_1) = ht(\overline{h_2})$

In this case there exists $n \geq 0$ such that $h_1^n = c^{-1} \circ h_3$, where $\overline{h_2} = h_3 \circ h_4$. Hence,
$ht(h_1) = ht(h_4 \circ h_3)$, both $h_1$ and $h_4 \circ h_3$ are cyclically reduced
and $c(h_1^k, (h_4 \circ h_3)^m) \geqslant |h_1| + |h_4 \circ h_3|$ for some $k, m > 0$.
By Lemma \ref{le:LS} we have $[h_1,h_4 \circ h_3] = \varepsilon$, that is, $C_G(h_1) = C_G(h_4 \circ h_3)$.
It follows
$$C_G(h_2) = c^{-1} \ast C_G(\overline{h_2}) \ast c = c^{-1} \ast C_G(h_3 \circ h_4) \ast c =$$
$$= c^{-1} \ast h_4^{-1} \ast C_G(h_4 \circ h_3) \ast h_4 \ast c = c^{-1} \ast h_4^{-1}
\ast C_G(h_1) \ast h_4 \ast c.$$
On the other hand, $h_1^n = c^{-1} \circ h_3$, so
$$C_G(h_2) = (c^{-1} \ast h_4^{-1} \ast h_3^{-1} \ast c) \ast C_G(h_1) \ast (c^{-1} \ast h_3
\ast h_4 \ast c)$$
$$= h_2^{-1} \ast C_G(h_1) \ast h_2,$$
that is, $C_G(h_1) = C_G(h_2)$ and $[h_1,h_2] = \varepsilon$.

\item[{\bf (b)}] $ht(h_1) < ht(\overline{h_2})$

There exists $n \geq 0$ such that $c^{-1} = h_1^n \circ h_3$, where $h_1 = h_3 \circ h_4$. Hence, from $ht(h_1) <
ht(\overline{h_2})$ it follows that
$\overline{h_2}$ contains any positive power of $h_4 \circ h_3$ as an initial subword.

Observe that $f^{-1} \ast h_1 \ast f = f_1^{-1} \ast (c \ast h_1 \circ c^{-1}) \ast f_1 =
f_1^{-1} \ast (c \ast (h_3 \circ h_4) \ast c^{-1}) \ast f_1 =
f_1^{-1} \ast (h_4 \circ h_3) \ast f_1$ and $ht(f_1^{-1} \ast (h_4 \circ h_3) \ast f_1) < ht(f_1)$.
From $ht(f_1^{-1} \ast (h_4 \circ h_3) \ast f_1) < ht(f_1)$ and $ht(f_1^{-1} \ast \overline{h_2} \ast f_1) < ht(f_1)$
we get $ht(f_1^{-1} \ast ((h_4 \circ h_3) \circ \overline{h_2}) \ast f_1) < ht(f_1)$.
It follows that $f_1$ contains any positive power of $(h_4 \circ h_3) \circ \overline{h_2}$
as an initial subword and $ht((h_4 \circ h_3) \circ \overline{h_2}) = ht(\overline{h_2})$.
By Lemma \ref{le:LS} we have $[(h_4 \circ h_3) \circ \overline{h_2},\overline{h_2}] = \varepsilon$ or,
equivalently, $[h_4 \circ h_3,\overline{h_2}] = \varepsilon$. That is, $C_G(\overline{h_2}) = C_G(h_4 \circ h_3)$.
Now, $h_4 \ast C_G(h_3 \circ h_4) \ast h_4^{-1} = C_G(\overline{h_2})$ and
$C_G(h_1) = h_4^{-1} \ast C_G(\overline{h_2}) \ast h_4$. From $c^{-1} \circ h_4 = h_1^{n+1}$
it follows that $h_4 = c \ast h_1^{n+1}$ and
$$C_G(h_1) = h_1^{-(n+1)} \ast c^{-1} \ast C_G(\overline{h_2}) \ast c \ast h_1^{(n+1)} =
h_1^{-(n+1)} \ast C_G(h_2) \ast h_1^{(n+1)}.$$
Hence, $C_G(h_1) = C_G(h_2)$ and $[h_1,h_2] = \varepsilon$.

\item[{\bf (c)}] $ht(h_1) > ht(\overline{h_2})$

There exists $n \geq 0$ such that $c^{-1} = h_1^n \circ h_3$, where $h_1 = h_3 \circ h_4$ and $h_4 \circ h_3$
contains any positive power of $\overline{h_2}$ as an initial segment. Similarly to
{\bf (b)} we get $[h_4 \circ h_3,\overline{h_2}] = \varepsilon$ and derive $C_G(h_1) = C_G(h_2)$.

\end{enumerate}

\item[{\bf (2)}] $ht(c) > ht(h_1)$

It follows that $c^{-1}$ contains any positive power of $h_1$ as an initial segment.
Consider $f^{-1} \ast (h_1 \circ h_2) \ast f$. We have
$$f^{-1} \ast (h_1 \circ h_2) \ast f = f^{-1} \ast (h_1 \circ c^{-1} \circ \overline{h_2}
\circ c) \ast f$$
$$ = (f_1^{-1} \circ c) \ast (h_1 \circ c^{-1} \circ \overline{h_2} \circ c) \ast
(c^{-1} \circ f_1) = (f_1^{-1} \circ c) \ast (h_1 \circ c^{-1} \circ \overline{h_2} \circ f_1).$$
Hence, $h_1 \circ c^{-1}$ cancels completely in $(f_1^{-1} \circ c) \ast (h_1 \circ c^{-1})$ and
$c \ast (h_1 \circ c^{-1}) = h_3$ so that $|h_1| = |h_3|$. Thus,
$(f_1^{-1} \circ c) \ast (h_1 \circ c^{-1} \circ \overline{h_2} \circ f_1) = f_1^{-1} \ast (h_3 \circ
\overline{h_2} \circ f_1)$ and since $ht(f_1^{-1} \ast (h_3 \circ \overline{h_2}) \ast f_1) < ht(f_1)$,
by Lemma \ref{le:0}, $f_1$ contains any positive power of $h_3 \circ \overline{h_2}$
as an initial subword. On the other hand,
$ht(f^{-1} \ast h_1 \ast f) < ht(f)$, so, $ht((f_1^{-1} \circ c) \ast h_1 \ast (c^{-1} \circ f_1)) =
ht(f_1^{-1} \ast (h_3 \circ f_1)) < ht(f_1)$, and it follows that $f_1$ contains any positive power of $h_3$
as an initial segment. Using Lemma \ref{le:LS} we get $[h_3,\overline{h_2}] = \varepsilon$ and
$C_G(\overline{h_2}) = C_G(h_3) = c \ast C_G(h_1) \ast c^{-1}$. Hence, $C_G(h_2) = c^{-1} \ast
C_G(\overline{h_2}) \ast c = C_G(h_1)$ and $[h_1,h_2] = \varepsilon$.
\end{enumerate}

\end{proof}

\begin{lemma}
\label{le:2} Let $f,h_1 {\not = \epsilon} \in G$ be such that $f$ is cyclically reduced and $ht(h_1) < ht(f)$.
If $ht(f^{-1} \ast h_1 \ast f) < ht(f)$ and $[h_1,h_2] = \varepsilon$, where $ht(h_2) < ht(f)$, then $ht(f^{-1} \ast h_2 \ast f)
< ht(f)$.
\end{lemma}
\begin{proof} Without loss of generality we can assume $h_1$ to be cyclically reduced. Indeed,
if $h_1 = c^{-1} \circ \overline{h_1} \circ c$ is the cyclic decomposition of $h_1$ then
by Lemma \ref{le:cycl}, $h_2 = c^{-1} \circ \overline{h_2} \circ c$ is the cyclic decomposition of $h_2$ and from
$ht(f^{-1} \ast (c^{-1} \circ \overline{h_1} \circ c) \ast f) < ht(f)$ it follows that
$f = c^{-1} \circ f_1$. Hence, $f_2 = f_1 \circ c^{-1}$ is cyclically reduced and from
$f_2^{-1} \ast \overline{h_1} \ast f_2 = (c \circ f_1^{-1}) \ast \overline{h_1}
\ast (f_1 \circ c^{-1}) = c \ast (f_1^{-1} \ast \overline{h_1} \ast f_1) \ast c^{-1}$,
we get $ht(f_2^{-1} \ast \overline{h_1} \ast f_2) < ht(f_2)$,
that is, we can consider the triple $f_2, \overline{h_1}, \overline{h_2}$ instead of $f, h_1, h_2$.

Since $h_1$ is cyclically reduced it follows that in the product $f^{-1} \ast h_1 \ast f$
either $f^{-1} \ast h_1 = f^{-1} \circ h_1$ or $h_1 \ast f = h_1 \circ f$. Assume the latter.
Now, from $ht(f^{-1} \ast h_1 \ast f) < ht(f)$, by Lemma \ref{le:0}, $f$ contains any
positive power of $h_1$ as an initial segment.

Since, $[h_1,h_2] = \varepsilon$ then by Lemma \ref{le:cycl}, $h_2$ is cyclically reduced and we have
$$(f^{-1} \ast h_1 \ast f) \ast (f^{-1} \ast h_2 \ast f) = (f^{-1} \ast h_2 \ast f) \ast (f^{-1} \ast h_1 \ast f).$$
Denote $g = f^{-1} \ast h_1 \ast f$ and let $g = c^{-1} \circ \overline{g} \circ c$
be the cyclic decomposition of $g$. Notice that $ht(c) < ht(f)$. Now,
$$(f^{-1} \ast h_2^{-1} \ast f) \ast g \ast (f^{-1} \ast h_2 \ast f) = g$$
and conjugating both sides by $c^{-1}$ we get
$$(f_1^{-1} \ast h_2^{-1} \ast f_1) \ast \overline{g} \ast (f_1^{-1} \ast h_2 \ast f_1) = \overline{g},$$
where $f = f_1 \circ c$.

Suppose the statement of the lemma does not hold, that is, $ht(f_1^{-1} \ast h_2 \ast f_1) = ht(f)$.
Then $ht(f_1^{-1} \ast h_2^{-1} \ast f_1) = ht(f)$ and $f_1^{-1} \ast h_2^{-1} \ast
f_1 = f_2^{-1} \circ h_3^{-1} \circ f_2$, where $f_1 = d \circ f_2,\ ht(f_2) = ht(f),\ ht(h_3) < ht(f)$.
Since $[f_1^{-1} \ast h_2^{-1} \ast f_1, \overline{g}] = [f_2^{-1} \circ h_3^{-1} \circ f_2, \overline{g}] =
\varepsilon$, by Lemma \ref{le:cycl} we have $f_2 = \epsilon$. Then $ht(f_1^{-1} \ast h_2 \ast f_1) = ht(h_3) < ht(f)$,
which contradicts the assumption.

\end{proof}

\begin{prop}
\label{pr:1} \cite{AB}
If $G < CDR(A,X)$ then for any $g \in G$, its centralizer $C_G(g)$ is a subgroup of
$A$. In particular, if $A = \mathbb{Z}^n$ then $C_G(g)$ is a free abelian group of
rank not more than $n$.
\end{prop}

\section{Complete $\mathbb{Z}^n$-free groups}
\label{sec:main}

In this section we fix a finitely generated group $G$ which has a
free regular length function with values in $\mathbb{Z}^n,\ n \in
\mathbb{Z}$ (with the right lexicographic order). Due to Theorem
\ref{chis-cor-1} we may and will  view $G$ as a subgroup of
$CDR(\mathbb{Z}^n,X)$ for an appropriate set $X$. Therefore, elements of $G$
are infinite words from $CDR(\mathbb{Z}^n,X)$, multiplication in $G$ is the
multiplication $\ast$ of infinite words, and the regular length function is the
standard length $| \cdot |$ of infinite words.

In the ordered group $\mathbb{Z}^n$ with basis $a_1, \ldots, a_n$ the subgroups
$E_k = \langle a_1, \ldots,$ $a_k\rangle$ are convex, and every non-trivial convex
subgroup is equal to $E_k$ for some $k$, so
$$0 = E_0 < E_1 < \cdots < E_n,$$
is the complete chain of convex subgroups of $\mathbb{Z}^n$. Recall,
that the height $ht(g)$ of a word $g \in G$ is equal to $k$ if
$|g| \in E_k - E_{k-1}$ (see \cite{KMS}). Since
$|g\ast h| \leq |g| + |h|$ and $|g^{-1}| = |g|$ one has for any $f,g \in G$:
\begin{enumerate}
\item $ht(f \ast g) \leqslant \max\{ht(f),ht(g)\}$,
\item $ht(g) = ht(g^{-1})$.
\end{enumerate}
We will assume that there is an element $g \in G$ with $ht(g) = n$,
otherwise, the length function on $G$ has  values in $\MZ^{n-1}$, in which
case we  replace $\MZ^n$ with $\MZ^{n-1}$.
For any $k \in [1,n]$
$$G_k = \{g \in G \mid ht(g) \leqslant k\}$$
is a subgroup of $G$ and
$$1 = G_0 < G_1 < \cdots < G_n = G.$$
Observe, that if $l:G \to \Lambda$ is a Lyndon length function with values
in some ordered abelian group $\Lambda$ and $\mu:\Lambda \to \Lambda^\prime$
is a homomorphism of ordered abelian groups then the composition
$l^\prime = \mu \circ l$ gives a Lyndon length function
$l^\prime : G \to \Lambda ^\prime$. In particular, since $E_{n-1}$ is a
convex subgroup of $\MZ^n$ then the canonical projection $\pi_n:\MZ^n \to \MZ$
such that  $\pi_n(x_1,x_2,\ldots,x_n) = x_n$ is an ordered homomorphism, so
the composition $\pi_n \circ |\cdot|$ gives a Lyndon length function
$\lambda: G \to \MZ$ such that $\lambda(g) = \pi_n(|g|)$. Notice also that
if $u = g \circ h$ then $\lambda(u) = \lambda(g) + \lambda(h)$ for any
$g, h, u \in G$.

Our goal in this section is to describe the algebraic structure  of $G$ in
terms of HNN-extensions of a very particular type. To this end in
Subsection \ref{subsec:elem-moves} we introduce Nielsen-like moves and
show that every finite set $Y$ of generators of $G$ can be transformed
by a sequence of elementary moves into a set $Z$ with a minimal total
$\lambda$-length. Then, in Subsection \ref{subsec:minimal}, we show that any
minimal finite set of generators $Z$ of $G$ gives an HNN-splitting of $G$.
Here we use the technique of Stallings' pregroups and their universal groups.
Finally, in Subsection \ref{subsec:structure} we use the HNN-splittings obtained
earlier to prove the main structure theorem on groups acting freely on
$\MZ^n$-trees.

\subsection{Elementary transformations of infinite words}
\label{subsec:elem-moves}

In this section we describe  an analog of Nielsen reduction  in the group $G$.
Since $G$ is complete the Nielsen reduced sets have much stronger non-cancellation
properties then usual and the transformations are simpler. On the other hand, since
$\MZ^n$ is non-Archimedean (for $n>1$) the reduction process is more cumbersome, it
goes in stages along the complete series of convex subgroups in $\MZ^n$.

For a finite subset  $Y$ of $G$ define its $\lambda$-length as
$$|Y|_\lambda = \sum_{g \in Y} \lambda(g).$$
If $Y$ is a generating set of $G$ then  $|Y|_\lambda >0$, otherwise $G = G_{n-1}$.
It follows that
$$Y = Y_+ \cup Y_0,$$
where
$$Y_+ = \{g \in Y \mid \lambda(g) > 0\},\ Y_{0} = \{g \in Y \mid \lambda(g) = 0\}.$$
Obviously, $|Y|_\lambda = |Y_+|_\lambda$ and $\langle Y_{0} \rangle$ is a
finitely generated subgroup of $G_{n-1}$.

Let $Y$ be a finite generating set for $G$. Assuming $Y = Y^{-1}$ we define three types of elementary transformations of $Y$
(below when we add/remove an element we also add/remove its inverse to keep the property $Y = Y^{-1}$).

\medskip {\bf Transformation $\mu$}. Let  $f,g \in Y_+,\ f \neq g$,
$h \in \langle Y_0 \rangle$,  $u = com(f,h\ast g)$, and $\lambda(u) >0$. Then
$f = u \circ w_1,\ h \ast g = u \circ w_2$ for some $u, w_1, w_2$ from $G$ (since $G$ is complete).
Put
$$\mu_{f,g,h}(Y) = (Y - \{f,g\}) \bigcup \{w_1, w_2, u\} \bigcup \{f^{-1} \ast h \ast g \mid
{\rm if}\ \lambda(f^{-1} \ast h \ast g) = 0\}$$
if $f \neq g^{-1}$, and
$$\mu_{f,g,h}(Y) = (Y - \{g\}) \bigcup \{u, w_2 \ast u\}$$
if $f = g^{-1}$.

\begin{lemma}
\label{le:mu}
In the notation above
\begin{enumerate}
\item[(1)] $\langle Y \rangle = \langle \mu_{f,g,h}(Y) \rangle$,
\item[(2)] $|\mu_{f,g,h}(Y)|_\lambda < |Y|_\lambda$.
\end{enumerate}
\end{lemma}
\begin{proof} Suppose $f \neq g^{-1}$. Then (1) is obvious. Next,
$$\lambda(f) + \lambda(g) = 2\lambda(u) + \lambda(w_1) + \lambda(w_2) > \lambda(u) +
\lambda(w_1) + \lambda(w_2),$$
that is, $|\mu_{f,g,h}(Y)|_\lambda < |Y|_\lambda$ and (2) follows.

Now, suppose $f = g^{-1}$. Then (1) is obvious again: $\langle \mu_{f,g,h}(Y) \rangle$ contains
$u,\ w_2 \ast u = u^{-1} \ast (h \ast g) \ast u$ and $h$, hence it contains $g$ as well. To see (2) let
$h = h' \circ h'',\ g = (h'')^{-1} \circ g'$, so that $h \ast g = h' \circ g'$. We have
$$h' \circ g' = u \circ w_2,\ g^{-1} = ((h'')^{-1} \circ g')^{-1} = u \circ w_1,$$
where $u = com(g^{-1},h \ast g)$ and $w_1, w_2 \in G$. That is,
$$h' \circ g' = u \circ w_2,\ (h'')^{-1} \circ g' = w_1^{-1} \circ u^{-1}.$$
First of all notice that $\lambda(g) = \lambda(g') = \lambda(u) + \lambda(w_1) = \lambda(u) + \lambda(w_2)$,
so, since $\lambda(u) > 0$ it follows that $\lambda(w_1) = \lambda(w_2) < \lambda(g')$, hence, in particular,
$|w_2| < |g'|$.

Consider the terminal subwords of length $|g'|$ in both $u \circ w_2$ and $w_1^{-1} \circ u^{-1}$.
Since they are the same then $|u| < |w_2|$. Indeed, if $|u| > |w_2|$ then some terminal subword
of $u$ is equal to its inverse, and if $|u| = |w_2|$ then $w_2 = u^{-1}$ and $h' \circ g' = u \circ u^{-1}$,
so, both cases are impossible. Hence, $w_2$ has $u^{-1}$ as a terminal
subword, that is, $w_2 = w \circ u^{-1}$ and
$$h' \circ g' = u \circ w_2 = u \circ w \circ u^{-1}$$
where
$$\lambda(g) = \lambda(g') = 2\lambda(u) + \lambda(w) > \lambda(u) + \lambda(w) = \lambda(u) + \lambda(u^{-1} \ast (h \ast g) \ast u).$$
Thus, $|\mu_{f,g,h}(Y)|_\lambda < |Y|_\lambda$.

\end{proof}

\medskip

{\bf Transformation $\eta$}. Let $f \in Y_+$ be
such that $\lambda(f) > \lambda(com(f,h\ast f)) > 0$ for some $h \in \langle Y_0
\rangle$. Then $f = u \circ f_1,\ h \ast f = u \circ f_2,\ \lambda(u) > 0$.
Define
$$\eta_{f,h}(Y) = (Y - \{f\}) \bigcup \{f_1, u, u^{-1} \ast h \ast u\}.$$
Notice, that $f = (h^{-1} \ast u) \circ f_2 = u \circ f_1$ hence $\lambda(f_2) = \lambda(f_1) > 0$.
On the other hand $f_2 = (u^{-1}\ast h\ast u) \ast f_1$ and it follows that
$\lambda(u^{-1} \ast h \ast u) = 0$.

\begin{lemma}
\label{le:eta}
In the notation above
\begin{enumerate}
\item[(1)] $\langle Y \rangle = \langle \eta_{f,h}(Y) \rangle$,
\item[(2)] if $f_1 \neq u$ then $|\eta_{f,h}(Y)|_\lambda = |Y|_\lambda$ and $|{\eta_{f,h}(Y)}_+| > |Y_+|$,
\item[(3)] if $f_1 = u$ then $|\eta_{f,h}(Y)|_\lambda < |Y|_\lambda$ and $|{\eta_{f,h}(Y)}_+| = |Y_+|$.
\end{enumerate}
\end{lemma}
\begin{proof} Directly from the construction.
\end{proof}

\medskip

{\bf Transformation $\nu$}. Let $f \in
Y_+$ be not cyclically reduced. Then $f = c^{-1} \circ \overline{f} \circ c$,
where $c \neq 1$ and $\overline{f}$ is cyclically reduced. In this case
$c^{-1} = com(f,f^{-1})$, hence (since $G$ is complete) $c, \overline{f} \in G$.
Put
$$\nu_f(Y) = (Y - \{f\}) \bigcup \{c, \overline{f}\}.$$

\begin{lemma}
\label{le:nu}
In the notation above
\begin{enumerate}
\item[(1)] $\langle Y \rangle = \langle \nu_f(Y) \rangle$,
\item[(2)] $|\nu_f(Y)|_\lambda \leqslant |Y|_\lambda$,
\item[(3)] $|{\nu_f(Y)}_+| \geq |Y_+|$
\end{enumerate}
\end{lemma}
\begin{proof} Obvious.
\end{proof}

We write $Y \to Y'$ ($Y \to^\ast Y'$) if $Y'$ is obtained from $Y$
by a single (finitely many) elementary transformation, that is, $\to^\ast$
is the transitive closure of the relation $\to$. We call a generating set
$Y$ of $G$ {\it transformation-reduced}  if none of the transformations
$\mu, \eta, \nu$ can be applied to $Y$. Recall that the binary relation
$\to^\ast$ is called {\it terminating} if there is no an infinite sequence
of finite subsets $Y_i, i\in \MN$, of $G$ such that $Y_i \to Y_{i+1}$ for
every $i \in \MN$, i.e., every rewriting system $Y_1 \to Y_2 \to \ldots $
is finite. We say that $\to^\ast$ is {\em uniformly terminating} if for
every finite set $Y$ of $G$ there is a  natural number $n_Y$ such that every
rewriting system starting at $Y$ terminates in at most $n_Y$ steps.

\begin{prop}
\label{pr:min}
The following hold:
\begin{enumerate}
\item [1)] The relation $\to^\ast$ is uniformly terminating. Moreover, for any
finite subset $Y$ of $G$ one has $n_Y \leq (|Y|_\lambda)^2$.
\item [(2)] If $Z$ is  a transformation-reduced finite subset of $G$ then:
\begin{itemize}
\item[(a)] all elements of $Z_+$ are cyclically reduced;
\item[(b)] if $f,g \in Z_+^{\pm 1},\ f \neq g$ then $\lambda(com(f,h\ast g)) = 0$
for any $h \in \langle Z_0 \rangle$;
\item[(c)] if $f \in Z_+^{\pm 1}$ and $\lambda(com(f,h\ast f)) > 0$ for some
$h \in \langle Z_0 \rangle$ then $\lambda(com(f,h\ast f)) = \lambda(f)$.
\end{itemize}
\item [(3)] If $Z$  is a transformation-reduced finite subset of $G$ then one can add to $Z$
finitely many elements $h_1, \ldots, h_m \in G_{n-1}$ such that $T = Z \cup \{h_1, \ldots, h_m\}$
is transformation-reduced and satisfies the following condition
\begin{itemize}
\item[(d)] if $f \in T_+^{\pm 1}$ and $\lambda(com(f,h\ast f)) > 0$ for some
$h \in \langle T_0 \rangle$ then $\lambda(com(f,h\ast f)) = \lambda(f)$ and
$f^{-1} \ast h \ast f \in \langle T_0 \rangle$.
\end{itemize}
\end{enumerate}
\end{prop}
\begin{proof} The existence of a finite transformation-reduced $Z$ which is obtained from $Y$ by
finitely many applications of $\mu, \eta, \nu$ and which satisfies (a), (b)
and (c) follows from Lemma \ref{le:mu}, \ref{le:eta} and \ref{le:nu},
decreasing step by step $|Y|_\lambda$ while increasing $|Y|$ which
is bounded from above by $(|Y|_\lambda)^2$. Finally, observe that if $f \in Z_+^{\pm 1}$
and $\lambda(com(f,h\ast f)) > 0$ for some $h \in \langle Z_0 \rangle$ then
$\lambda(f^{-1} \ast h \ast f) = 0$ by (c), that is, $f^{-1} \ast h \ast f \in G_{n-1}$.
Now, from Lemma \ref{le:1} and Lemma \ref{le:2}, $f^{-1} \ast h \ast f \in G_{n-1}$ if and only if
$g \in C_{G_{n-1}}(h)$, and in order for $Z$ to satisfy (d), it is enough to add
$C_{G_{n-1}}(h)$ and $f^{-1} \ast C_{G_{n-1}}(h) \ast f$ to $Z_0$. By Proposition
\ref{pr:1} both centralizers are finitely generated and it is enough to add only
finitely many elements to $Z_0$.

\end{proof}

\begin{defn}
A finite set $Y$ of $G$ is called reduced if it satisfies the conditions (a) - (d) from Proposition \ref{pr:min}.
\end{defn}

\subsection{Minimal sets of generators and pregroups}
\label{subsec:minimal}

Let $Z$ be a finite reduced generating set of $G$.
Put
$$P_Z = \{g \ast f \ast h \mid f \in Z_+^{\pm 1},\ g,h \in \langle Z_0
\rangle \} \cup \langle Z_0 \rangle.$$
Multiplication $\ast$ induces a partial multiplication (which we again
denote by $\ast$) on $P_Z$ so that for $p, q \in P_Z$ the product
$p \ast q$ is defined in $P_Z$ if and only if $p \ast q \in P_Z$.
Notice, that $P_Z$ is closed under inversion.

\begin{lemma}
\label{le:claim}
Let $x = h_1(x) \ast f_x \ast h_2(x),\ y = h_1(y) \ast
f_y \ast h_2(y) \in P_Z$, where $h_i(x), h_i(y) \in \langle Z_0
\rangle,\ i = 1,2$ and $f_x, f_y \in Z_+^{\pm 1}$. Then $x \ast y \in
P_Z$ if and only if $f_x = f_y^{-1}$ and $f_x \ast (h_2(x) \ast
h_1(y)) \ast f_x^{-1} \in \langle Z_0 \rangle$.
\end{lemma}
\begin{proof} If $f_x \ast (h_2(x) \ast h_1(y)) \ast f_y \in
\langle Z_0 \rangle$ then since $Z$ is reduced we have $f_x =
f_y^{-1}$ and obviously $x \ast y \in \langle Z_0 \rangle \subset P_Z$.

Now, assume $x \ast y \in P_Z$. It follows
$$x \ast y = (h_1(x) \ast f_x \ast h_2(x)) \ast (h_1(y) \ast f_y \ast
h_2(y)) = h_1(z) \ast f_z \ast h_2(z) \in P_Z,$$
or, in other words
$$(h_1(x) \ast f_x \ast h_2(x)) \ast (h_1(y) \ast f_y \ast h_2(y))
\ast (h_2(z)^{-1} \ast f_z^{-1} \ast h_1(z)^{-1} = \varepsilon.$$
Thus, either we have $\lambda(com(f_x^{-1},(h_2(x) \ast h_1(y)) \ast f_y)) > 0$
or $\lambda(com(f_y^{-1}, (h_2(y) \ast h_2(z)^{-1}) \ast f_z^{-1}))
> 0$. Suppose the former is true (the latter case is considered similarly).
Since $Z$ is reduced it follows $f_x^{-1} = f_y$ and $f_x \ast (h_2(x) \ast
h_1(y)) \ast f_x^{-1} \in \langle Z_0 \rangle$.

\end{proof}

Now we are ready to prove the main technical result of this section.

\begin{theorem}
\label{pr:pregroup}
Let $G$ be a finitely generated group with a free length function in $\MZ^n$. Then:
\begin{enumerate}
\item[(1)] $P_Z$ forms a pregroup with respect to the multiplication $\ast$ and inversion;
\item[(2)] the inclusion $P_Z \to G$ extends to the group isomorphism $U(P_Z) \to G$,
where $U(P)$ is the universal group of $P_Z$;
\item[(3)] if $(g_1,\ldots,g_k)$ is a reduced $P_Z$-sequence for an
element $g \in G$ then
$$\lambda(g) = \sum^k_{i=1} \lambda(g_i).$$
\end{enumerate}
\end{theorem}

\begin{proof} Observe that $P_Z = P_Z^{-1} \subset G$ generates $G$ and every $g \in G$ corresponds
to a finite reduced $P_Z$-sequence
$$(u_1, u_2, \ldots, u_k),$$
where $u_i \in P_Z,\ i \in [1,k],\ u_i \ast u_{i+1} \notin P_Z,\ i \in [1,k-1]$
and $g = u_1 \ast u_2 \ast \cdots \ast u_k$ in $G$. By Theorem 2, \cite{R2}, to prove that $P_Z$ is a pregroup and the
inclusion $P_Z \to G$ extends to the isomorphism $U(P_Z) \to G$ it is enough to show
that all reduced $P_Z$-sequences representing the same element have the same $P_Z$-length.

Suppose two reduced $P_Z$-sequences
$$(u_1, u_2, \ldots, u_k),\ (v_1, v_2, \ldots, v_n)$$
represent the same element $g \in G$. That is,
$$(u_1 \ast \cdots \ast u_k) \ast (v_1 \ast \cdots \ast v_n)^{-1} = \varepsilon.$$
We use the induction on $k+n$ to show that $k = n$. If the $P_Z$-sequence
$$(u_1, \ldots, u_k, v_n^{-1}, \ldots, v_1^{-1})$$
is reduced then
$$u_1 \ast \ldots \ast u_k \ast v_n^{-1} \ast \ldots \ast v_1^{-1} \neq \varepsilon$$
because $Z$ is a reduced set. Hence,
$$(u_1, \ldots, u_k, v_n^{-1}, \ldots, v_1^{-1})$$
is not reduced and $u_k \ast v_n^{-1} \in P_Z$. If $u_k = h_1 \ast f_1 \ast g_1,\
v_n = h_2 \ast f_2 \ast g_2$, where $h_i, g_i \in \langle Z_0 \rangle$
and $f_i \in Z_+^{\pm 1},\ i = 1,2$ then by Lemma \ref{le:claim} $f_1 = f_2$ and
$f_1 \ast (g_1 \ast g_2^{-1}) \ast f_2^{-1} = c \in \langle Z_0 \rangle$. It follows
that
$$(u_1, u_2, \ldots, u_{k-1} \ast (h_1 \ast c \ast h_2^{-1})),\ (v_1, v_2, \ldots, v_{n-1})$$
represent the same element $g \ast v_n^{-1} \in G$ and the sum of their lengths is less
than $k+n$, so the result follows by induction. Hence, (1) and (2) follow.

\smallskip

Finally we prove (3).

If $g_i = h_1(g_i) \ast f_{g_i} \ast h_2(g_i),\ i \in [1,k]$ then $\lambda(g_i) = \lambda(f_{g_i})$ because
$\lambda(h_1(g_i)) = \lambda(h_2(g_i)) = 0$. On the other hand, since $Z$ is reduced and
$(g_1,\ldots,g_k)$ is a reduced $P_Z$-sequence then
$\lambda(com(g_i^{-1},g_{i+1})) = 0$ for $i \in [1,k-1]$. In other words $\lambda(g_i \ast g_{i+1}) = \lambda(g_i) +
\lambda(g_{i+1})$ and the result follows.

\end{proof}

\begin{corollary}
\label{co:1}
$G_{n-1} = \langle Z_0 \rangle$.
\end{corollary}
\begin{proof} From Theorem \ref{pr:pregroup} it follows that
$G = U(P_Z)$ and every element $g$ of $G$ can be
represented as a reduced $P_Z$-sequence $g = (g_1,\ldots,g_k)$,
where $g_i \in P_Z$ and $g_i \ast g_{i+1} \notin P_Z$ for any
$i \in [1,k-1]$. It follows that $g_i \notin \langle Z_0
\rangle$ for any $i \in [1,k]$. Hence, if $g \in G_{n-1}$ then
$\lambda(g) = 0$ and it can only be represented by a reduced
$P_Z$-sequence of length $1$. So, the statement follows from
Theorem \ref{pr:pregroup}, (3).

\end{proof}

\subsection{Algebraic structure of complete $\mathbb{Z}^n$-free
groups}
\label{subsec:structure}

\begin{theorem}
\label{th:main}
Let $G$ be a finitely generated group with a regular free Lyndon
length function in $\mathbb{Z}^n$ and let $Z$ be a reduced
generating set for $G$. Then $G$ has the following
presentation
$$G = \langle H, Y \mid t_i^{-1} C_H(u_{t_i}) t_i =
C_H(v_{t_i}),\ t_i \in Y^{\pm 1} \rangle,$$
where $Y = Z_+$ is finite, $H = G_{n-1} = \langle Z_0 \rangle$ is
finitely generated and $C_H(u_{t_i})$ $C_H(v_{t_i})$
are either trivial or finitely generated free abelian subgroups of
$H$. Moreover, $H$ has a regular free Lyndon length function in
$\mathbb{Z}^{n-1}$.
\end{theorem}
\begin{proof} From Theorem \ref{pr:pregroup} it follows that
$G = U(P_Z)$, where
$$P_Z = \{g \ast f \ast h \mid f \in Z_+^{\pm 1},\ g,h \in
\langle Z_0 \rangle \} \cup \langle Z_0 \rangle.$$
It follows that every element $g$ of $G$ can be represented as a
reduced $P_Z$-sequence $g = (g_1,\ldots,g_k)$, where $g_i \in P_Z$
and $g_i \ast g_{i+1} \notin P_Z$ for any $i \in [1,k-1]$. It
follows that $g_i \notin \langle Z_0 \rangle$ for any $i \in
[1,k]$. In fact, we have
$$G = U(G) = \langle P_Z \mid x y = z,\ (x,y,z \in
P_Z\ {\rm and}\ x \ast y = z) \rangle.$$
Denote $H = \langle Z_0 \rangle$ and $Y = Z_+$. By Corollary
\ref{co:1} we have $H = G_{n-1}$.

At first observe that $P_Z$ is infinite but for each $p \in P_Z$
either $p \in H$ or $p = h_1(p) \ast f_p \ast h_2(p)$,
where $f_p \in Y^{\pm 1}$ and $h_i(p) \in H,\ i = 1,2$. Hence,
every $p \in P_Z$ can rewritten in terms of $Y^{\pm}$ and
finitely many generators of $H$. On the other hand, if
$x,y,z \in P_Z$ and $x \ast y = z$ then one of the three of them
is in $H$ and without loss of generality we can assume $z \in H$.
Hence, either $x,y$ are in $H$ too, or $x,y \notin H$ and
assuming $x = h_1(x) \ast f_x \ast h_2(x),\ y = h_1(y) \ast
f_y \ast h_2(y)$, where $h_i(x), h_i(y) \in H,\ i = 1,2,\ f_x,
f_y \in Y^{\pm 1}$ by Lemma \ref{le:claim}
we get $f_x = f_y^{-1}$ and $f_x \ast (h_2(x) \ast
h_1(y)) \ast f_x^{-1} \in H$. Hence, every relator $x y = z$, where
$x,y,z \in P_Z$ and $x \ast y = z$ can be rewritten as
$$f_x \ast u_{x,y} \ast f_x^{-1} = v_{x,y,z},$$
where $f_x \in Y^{\pm 1}$ and $u_{x,y}, v_{x,y,z} \in H$. By Lemma
\ref{le:2} it follows that for each $y \in Y^{\pm 1}$ there exists
$u_y \in H$ such that $y \ast u \ast y^{-1} \in H$ for $u \in H$ if
and only if $u \in C_H(u_y)$. Since $Z$ is reduced it follows that
for each $y \in Y^{\pm 1}$ both $C_H(u_y)$ and $y \ast C_H(u_y) \ast y^{-1}$
are in $H$ and also note that $y \ast C_H(u_y) \ast y^{-1}$ is a centralizer
of some element in $H$. Hence, every
$$f_x \ast u_{x,y} \ast f_x^{-1} = v_{x,y,z},$$
is a consequence of
$$f_x \ast u_x \ast f_x^{-1} = v_x,$$
where $u_x,v_x$ depend only on $f_x$. Thus,
$$G = \langle Y, H \mid t_i^{-1} C_H(u_{t_i}) t_i =
C_H(v_{t_i}),\ t_i \in Y^{\pm 1} \rangle,$$
where $Y$ is finite, $H$ is finitely generated and $C_H(u_y), C_H(v_y)$
are finitely generated abelian (see Proposition \ref{pr:1}).

Finally, we have to show that $H$ has a regular free Lyndon length
function in $\mathbb{Z}^{n-1}$. Indeed, since $H = G_{n-1} < G$ then the
free Lyndon length function with values in $\mathbb{Z}^{n-1}$ is automatically
induced on $H$. We just have to check if it is regular.

Take $g,h \in H$ and consider $com(g,h)$. Since the length function
on $G$ is regular then $com(g,h) \in G = U(P_Z)$ and $com(g,h)$
can be represented by the reduced $P_Z$-sequence $(g_1,\ldots,g_k)$.
By Theorem \ref{pr:pregroup}, (3) it follows that
$$\lambda(com(g,h)) = \sum^k_{i=1} \lambda(g_i).$$
But if $\lambda(com(g,h)) > 0$ then $\lambda(g),\lambda(h) > 0$ -
contradiction with the choice of $g$ and $h$. Hence, $\lambda(g_i) = 0,\
i \in [1,k]$ and it follows that $k = 1$. Thus, $com(g,h) = g_1 \in H$.
This completes the proof of the theorem.

\end{proof}

\begin{theorem}
\label{th:main2}
Let $G$ be a finitely generated group with a regular free Lyndon
length function in $\mathbb{Z}^n$. Then $G$ can be represented as
a union of a finite series of groups
$$G_1 < G_2 < \cdots < G_n = G,$$
where $G_1$ is a free group of finite rank, and
$$G_{i+1} = \langle G_i, s_{i,1},\ \ldots,\ s_{i,k_i} \mid s_{i,j}^{-1}\ C_{i,j}\ s_{i,j} = \phi_{i,j}(C_{i,j}) \rangle,$$
where for each $j \in [1,k_i],\ C_{i,j}$ and $\phi_{i,j}(C_{i,j})$ are cyclically reduced centralizers of $G_i$,
$\phi_{i,j}$ is an isomorphism, and the following conditions are satisfied:
\begin{enumerate}
\item[(1)] $C_{i,j} = \langle c^{(i,j)}_1, \ldots, c^{(i,j)}_{m_{i,j}} \rangle,\
\phi_{i,j}(C_{i,j}) = \langle d^{(i,j)}_1, \ldots, d^{(i,j)}_{m_{i,j}} \rangle$,
where $\phi_{i,j}(c^{(i,j)}_k) = d^{(i,j)}_k,\ k \in [1,m_{i,j}]$ and
$$ht(c^{(i,j)}_k) = ht(d^{(i,j)}_k) < ht(d^{(i,j)}_{k+1}) = ht(c^{(i,j)}_{k+1}),\
k \in [1,m_{i,j}-1],$$
$$ht(s_{i,j}) > ht(c^{(i,j)}_k),$$
\item[(2)]  $|\phi_{i,j}(w)| = |w|$ for any $w \in C_{i,j}$,
\item[(3)] $w$ is not conjugate to $\phi_{i,j}(w)^{-1}$ in $G_i$ for any $w \in C_{i,j}$,
\item[(4)] if $A, B \in \{C_{i,1}, \phi_{i,1}(C_{i,1}), \ldots, C_{i,k_i},
\phi_{i,k_i}(C_{i,k_i})\}$ then either $A = B$, or $A$ and $B$ are not conjugate
in $G_i$,
\item[(5)] $C_{i,j}$ can appear in the list
$$\{ C_{i,k}, \phi_{i,k}(C_{i,k}) \mid k \neq j \}$$
not more than twice.
\end{enumerate}
\end{theorem}
\begin{proof} Existence of the series
$$G_1 < G_2 < \cdots < G_n = G,$$
where $G_{i+1},\ i \in [1,n-1]$ can be obtained from $G_i$ by finitely many
HNN-extensions in which associated subgroups are maximal abelian
of finite rank follows by induction applying Theorem \ref{th:main}.
Also, observe that $G_1$ has a free length function with values in $\mathbb{Z}$, hence,
by the result of Lyndon \cite{L} it follows that $G_1$ is a free group. Moreover,
$G_1$ is of finite rank by Theorem \ref{th:main}.

Now, consider $G_{i+1}$. By Theorem \ref{th:main} we can assume that
\begin{equation}
\label{eq:1}
G_{i+1} = \langle G_i, t_1, t_2, \ldots, t_p \mid
t_j^{-1} C_{G_i}(u_{t_j}) t_j = C_{G_i}(v_{t_j}) \rangle,
\end{equation}
where
\begin{enumerate}
\item[(a)] all $t_j$ are cyclically reduced,
\item[(b)] $G_i = \langle Y \rangle,\ ht(t_j) > ht(G_i)$ and
$$Y \cup \{t_1, t_2, \ldots, t_p\}$$
is a reduced generating set for $G_{i+1}$.
\end{enumerate}
In particular, $Y \cup \{t_1, t_2, \ldots, t_p\}$ is reduced, that is, it
has the properties listed in Proposition \ref{pr:min}.

At first, we can assume that all $C_{G_i}(u_{t_j}),\ C_{G_i}(v_{t_j})$ are cyclically reduced.
Indeed, if not then by Lemma \ref{le:cycl} we have $C_{G_i}(u_{t_j}) = c^{-1} \circ B \circ c$,
where $B$ is cyclically reduced, $c \in G_i$ by regularity of
the length function on $G_i$, and
$$(t_j^{-1} \ast c^{-1}) \ast B \ast (c \ast t_j) = C_{G_i}(v_{t_j}).$$
Thus, we can substitute $t_j$ by $c \ast t_j,\ C_{G_i}(u_{t_j})$ by $B$, and
the same can be done for $C_{G_i}(v_{t_j})$.

Observe that conjugation by $t_j$ induces an isomorphism between
$C_{G_i}(u_{t_j})$ and $C_{G_i}(v_{t_j})$, and since we can assume both
centralizers to be cyclically reduced then from
$$t_j^{-1} \ast C_{G_i}(u_{t_j}) \ast t_j = C_{G_i}(v_{t_j})$$
it follows that for $a \in C_{G_i}(u_{t_j}),\ b \in C_{G_i}(v_{t_j})$
if $t_j^{-1} \ast a \ast t_j = b$ then $|a| = |b|$.
In particular, if
$$C_{G_i}(u_{t_j}) = \langle c^{(i,j)}_1, \ldots, c^{(i,j)}_{m_{i,j}} \rangle,$$
where we can assume $ht(c^{(i,j)}_k) < ht(c^{(i,j)}_{k+1})$ for $k \in
[1,m_{i,j}-1]$, then all $d^{(i,j)}_k = t_j^{-1} \ast c^{(i,j)}_k \ast t_j$
generate $C_{G_i}(v_{t_j})$ and $|c^{(i,j)}_k| = |d^{(i,j)}_k|$. This proves
(1) and (2).

\smallskip

Suppose there exist $w_1 \in C_{G_i}(u_{t_j})$ and $g \in G_i$ such
that $g^{-1} \ast w_1 \ast g = w_2^{-1}$, where $w_2 = \phi_i(w_1)
\in C_{G_i}(v_{t_j})$. Observe that either $ht(g) \leq ht(w_1) = ht(w_2)$
and in this case $w_1$ is a cyclic permutation of $w_2^{-1}$, or $ht(g) > ht(w_1)$.
In the latter case, $g$ has any positive power of $w_1^\delta, \delta \in \{1,-1\}$
as an initial subword and any positive power of $w_2^{-\delta}$ as a terminal subword. Without loss of
generality we can assume $\delta = 1$. Hence, $t_j \ast g^{-1} = t_j \circ g^{-1}$ and
$$(t_j \circ g^{-1})^{-1} \ast w_1 \ast (t_j \circ g^{-1}) = w_1^{-1}.$$
Consider $h = com(t_j \circ g^{-1},(t_j \circ g^{-1})^{-1})$. Observe that
$h \in G_i$ and $|h^{-1} \ast w_1 \ast h| = |w_1|$ (because
$|(t_j \circ g^{-1})^{-1} \ast w_1 \ast (t_j \circ g^{-1})| = |w_1|$).
Thus, if $w_3 = h^{-1} \ast w_1 \ast h$ then $h$ ends with any positive power of $w_3$. We have
$t_j \circ g^{-1} = h^{-1} \circ f \circ h$, where $f$ is cyclically reduced.
But at the same time we have $f^{-1} \ast w_3 \ast f = w_3^{-1}$ and this
produces a contradiction. Indeed, if $ht(f) \leq ht(w_3)$ then $w_3$ is a
cyclic permutation of $w_3^{-1}$ which is impossible. On the other hand, if
$ht(f) > ht(w_3)$ then $f$ has any power of $w_3^\alpha, \alpha \in
\{1,-1\}$ as an initial subword, and any power of $w_3^{-\alpha}$ as a
terminal subword - a contradiction with the fact that $f$ is cyclically
reduced. This proves (3).

\smallskip

To prove (4), assume that two centralizers from the list
$$C_{G_i}(u_{t_1}), \ldots, C_{G_i}(u_{t_p})$$
are conjugate in $G_i$. Denote $C_1 = C_{G_i}(u_{t_1}),\ C_2 =
C_{G_i}(u_{t_2})$ and let $C_1 = h^{-1} \ast C_2 \ast h$ for some $h \in G_i$.
Hence, in (\ref{eq:1}) every entry of $C_1$ can be substituted by
$h^{-1} \ast C_2 \ast h$ and some of the elements $t_1, t_2, \ldots, t_p$
can be changed accordingly.

\smallskip

Finally, assume that there exists $t_k \neq t_j$ such that $t_k^{-1} \ast C_{G_i}(u_{t_j})
\ast t_k \leq G_i$. Suppose $c(t_j, t_k) > 0$ and denote $z = com(t_j,t_k)$.
Observe that $ht(z) > ht(C_{G_i}(u_{t_j}))$. If $ht(z) < ht(t_j)$ then
$z$ conjugates $C_{G_i}(u_{t_j})$ into a cyclically reduced centralizer $A$
of $G_i$ and $z$ has any positive power of some $a \in A,\ ht(a) = ht(A)$
as a terminal subword. But then $ht(z^{-1} \ast t_j) = ht(z^{-1} \ast t_k) = ht(t_j)$
and since both $z^{-1} \ast t_j$ and $z^{-1} \ast t_k$ conjugate $A$ into
a cyclically reduced centralizer of $G_i$ it follows that $z^{-1} \ast t_j$ and
$z^{-1} \ast t_k$ have $a^{\pm 1}$ as an initial subword.
If $z^{-1} \ast t_j$ has $a$ as an initial subword and $z^{-1} \ast t_k$ has
$a^{-1}$ as an initial subword then $z \ast (z^{-1} \ast t_k) \neq
z \circ (z^{-1} \ast t_k)$, and we have a contradiction. If both $z^{-1} \ast t_j$ and
$z^{-1} \ast t_k$ have $a$ as an initial subword then $z$ cannot be $com(t_j,t_k)$
and again we have a contradiction. Thus, $ht(z) = ht(t_j)$, but it is possible only
if $t_j = t_k$ since $Y \cup \{t_1, t_2, \ldots, t_p\}$ is a minimal generating set, and
again we a contradiction with our choice of $t_k$. It follows that $c(t_j, t_k) = 0$
and if $t_j$ begins with $c \in C_{G_i}(u_{t_j})$ then $t_k$ begins with $c^{-1}$.
It also follows that there can be only one $t_k \neq t_j$ such that
$t_k^{-1} \ast C_{G_i}(u_{t_j}) \ast t_k \leq G_i$.

This completes the proof of the theorem.

\end{proof}

\section{Regular free actions of HNN-extensions}
\label{sec:embedding}

Let $H$ be a group with a regular free length function with values in
$\mathbb{Z}^n$. The main goal of this section is to prove the following result.

\begin{theorem}
\label{main4}
Let $H$ be a group with a regular free length function in $\mathbb{Z}[t]$.
Let $A$ and $B$ be centralizers in $H$ whose elements are cyclically reduced
and such that there exists an isomorphism $\phi : A \rightarrow B$ with the
following properties
\begin{enumerate}
\item $a$ is not conjugate to $\phi(a)^{-1}$ in $H$ for any $a \in A$,
\item $|\phi(a)| = |a|$ for any $a \in A$.
\end{enumerate}
Then the group
\begin{equation}
\label{eq:G}
G = \langle H, z \mid z^{-1} A z = B \rangle,
\end{equation}
has a regular free length function in $\mathbb{Z}[t]$ which extends
the length function on $H$.
\end{theorem}

\subsection{Cyclically reduced centralizers and attached
elements}

From Theorem \ref{th:main2}, $H$ is union of the chain
$$F(X) = H_1 < H_2 < \cdots < H_n = H,$$
where
$$H_{i+1} = \langle H_i, s_{i,1},\ldots,\ s_{i,k_i} \mid s_{i,j}^{-1}
C_{i,j} s_{i,j} = D_{i,j} \rangle,$$
$C_{i,j},\ D_{i,j}$ are maximal abelian subgroups of $H_i$, and
$ht(s_{i,j}) > ht(H_i)$ for any $i \in [1,n-1],\ j \in [1,k_i]$.

Let $K$ be a cyclically reduced centralizer in $H$. It is easy to see that
either $ht(K) = ht(H) = n$, or $ht(K) < ht(H)$ and $K$ is conjugate in
$H$ to a centralizer from $H_{n-1}$.

For a cyclically reduced centralizer $K$ of $H$ we define
$${\cal C}(K) = \{C_{i,j},\ D_{i,j}\} \cap \{ {\rm cyclically\ reduced\
centralizers\ conjugate\ to}\ K\ {\rm in}\ H\}.$$
Obviously, if ${\cal C}(K)$ is empty then $K$ is an infinite cyclic.
Now, for $C \in {\cal C}(K)$ we call $w$ from the list $s_{i,j},\ i
\in [1,n-1],\ j \in [1,k_i]$ {\em attached to $C$} if $ht(w) > ht(C)$
and $ht(w^{-1} \ast C \ast w) = ht(C)$. Observe that by Theorem
\ref{th:main2}, $C$ can have at most two attached elements of the same
height, and if $w_1,w_2$ are attached to $C$ and $ht(w_1) = ht(w_2)$ then
$w_1^{-1} \ast w_2 = w_1^{-1} \circ w_2$.

\begin{lemma}
\label{le:attached}
Let $K$ be a cyclically reduced centralizer of $H$, and let
${\cal C}(K)$ be empty. Let $a$ be a generator of $K$ of maximal
height. Then there is no element in $H$ which has any positive exponent of
$a^{\pm 1}$ as an initial subword.
\end{lemma}
\begin{proof} Suppose there exists $g \in H$ which starts with any
positive power of $a^{\pm 1}$. Consider $w = com(g, a \ast g)$ and
observe that we have $ht(w^{-1} \ast K \ast w) = ht(K)$. Without
loss of generality we assume $ht(w) = n$. Since $w \in H_n - H_{n-1}$,
there exists the following representation of $w$
$$w = h_1 \ast t_1^{\epsilon_1} \ast h_2 \ast \cdots \ast
t_k^{\epsilon_k} \ast h_{k+1},$$
where $t_i \in \{s_{n-1,1},\ldots,s_{n-1,k_{n-1}}\},\ i \in [1,k],\
h_j \in H_{n-1},\ j \in [1,k+1]$. We can assume $\epsilon_1 = 1$.

Denote $C = h_1^{-1} \ast K \ast h_1$, so $ht(t_1^{-1} \ast C \ast
t_1) = ht(C)$. At the same time there exists $A \in \{C_{i,j},D_{i,j}\},\
i \in [1,n-1],\ j \in [1,k_i]$ such that $ht(t_1^{-1} \ast A \ast
t_1) = ht(A)$. Thus by Lemma \ref{le:1} it follows that $C =
h_1^{-1} \ast K \ast h_1 = A$. But this implies that ${\cal C}(K)$ is
not empty - a contradiction.

\end{proof}

Below we are going to distinguish attached elements in the following way.
Suppose $C \in {\cal C}(K)$, and let $w$ be an element attached to $C$.
If $c$ is a generator of $C$ of maximal height then we call $w$ {\it
left-attached to $C$ with respect to $c$} if $c^{-1} \ast w = c^{-1}
\circ w$, and we call $w$ {\it right-attached to $C$ with respect to
$c$} if $c \ast w = c \circ w$.

\begin{lemma}
\label{le:attached2}
Let $K$ be a cyclically reduced centralizer of $H$, and let
$C \in {\cal C}(K)$. Let $c$ be a generator of $C$ of maximal height.
If there exists a right(left)-attached to $C$ with respect to $c$ element,
then there exists $D \in {\cal C}(K)$ and its generator $d$ of maximal height
such that $c$ is conjugate to $d$ in $H$ and $D$ does not have
right(left)-attached with respect to $d$ elements.
\end{lemma}
\begin{proof} Assume that there exists an element $t_1$ which is right-attached to $C$ with respect to $c$
(the case of a left-attached element is symmetric). Denote $D_1 = t_1^{-1}
\ast C \ast t_1 \in {\cal C}(K)$. Observe that $t_1^{-1} \ast (c \circ t_1)
= d_1 \in D_1$ and $ht(d_1) = ht(D_1) = ht(C)$. Also, it follows that $d_1$ is
not a proper power because otherwise $c$ is a proper power and this is
not possible. Hence, $d_1$ is a primitive element of $D_1$,
so, it is a generator of $D_1$ of maximal height. Now, observe that $t_1^{-1}$ is
a left-attached to $D_1$ with respect to $d_1$ element. If other attached to
$D_1$ elements are all left-attached with respect to $d_1$ then we are done.
Otherwise there exists a right-attached to $D_1$ with respect to $d_1$ element
$t_2$. Observe that $t_1 \ast t_2 = t_1 \circ t_2$. Denote $D_2 = t_2^{-1}
\ast (D_1 \circ t_2) \in {\cal C}(K)$ and repeat the argument above.

\smallskip

Since ${\cal C}(K)$ is finite then after a finite number of steps either we
find $D_p \in {\cal C}(K),\ D_p \neq C$ and its generator $d_p$ of maximal
height such that $D_p$ has only left-attached elements with respect to $d_p$,
or, we get $D_p = D_{p+m}$. Let us consider the latter case.
Without loss of generality we can assume $p = 1$. Thus, if we denote
$w = t_1 \circ t_2 \circ \cdots \circ t_m$ then $w^{-1} \ast C \ast w = C$,
that is, $w \in C$. But this gives us a contradiction because
$ht(w) > ht(C)$.

\end{proof}

\begin{lemma}
\label{le:attached3}
Let $K$ be a cyclically reduced centralizer of $H$, and let
$C \in {\cal C}(K)$. Let $c$ be a generator of $C$ of maximal height.
If there exists no right(left)-attached to $C$ with respect to $c$ element,
then there is no element $g \in H$ which has any positive exponent of $c$
as an initial(terminal) subword.
\end{lemma}
\begin{proof} Observe that it is enough to consider only the case when
there exists no right-attached to $C$ with respect to $c$ element (the
other case is symmetric).

The proof follows from the following claim.

\smallskip

{\bf Claim.} If there exists $g \in H$ which starts with any positive
exponent of $c$ then there exists $f \in H$ with the same
property and such that $ht(f) < ht(g)$.

\smallskip

Suppose there exists $g \in H$ which starts with any positive exponent of $c$.
Without loss of generality we assume $ht(g) = n$. Consider $w = com(g, c
\circ g)$. If $ht(w) < ht(g)$ then we are done. So, assume $w \in H_n -
H_{n-1}$.

Observe that we have $ht(w^{-1} \ast C \ast w) = ht(C)$ and without
loss of generality we can assume $\lambda(w)$ to be minimal
possible. Since $w \in H_n - H_{n-1}$, there exists
the following representation of $w$
$$w = h_1 \ast t_1^{\epsilon_1} \ast h_2 \ast \cdots \ast
t_k^{\epsilon_k} \ast h_{k+1},$$
where $t_i \in \{s_{n-1,1},\ldots,s_{n-1,k_{n-1}}\},\ i \in [1,k],\
h_j \in H_{n-1},\ j \in [1,k+1]$, and such that
$$\lambda(w) = \sum_{i=1}^k \lambda(t_i).$$
We can assume $\epsilon_1 = 1$. Then $ht((h_1 \ast t_1 \ast h_2)^{-1}
\ast C \ast (h_1 \ast t_1 \ast h_2)) = ht(C)$ and $\lambda(h_1 \ast t_1
\ast h_2) \leq \lambda(w)$. From our assumption it follows $\lambda(h_1
\ast t_1 \ast h_2) = \lambda(w)$ and $w = h_1 \ast t_1 \ast h_2$.

Denote $D = h_1^{-1} \ast C \ast h_1$, so $ht(t_1^{-1} \ast D \ast
t_1) = ht(D)$. At the same time there exists $A \in \{C_{i,j},D_{i,j}\},\
i \in [1,n-1],\ j \in [1,k_i]$ such that $ht(t_1^{-1} \ast A \ast
t_1) = ht(A)$. Thus by Lemma \ref{le:1} it follows that $D =
h_1^{-1} \ast C \ast h_1 = A$ and in particular $ht(C) = ht(A)$.

If $ht(h_1) \leq ht(C)$ then from Theorem \ref{th:main2} it
follows that $A = C$ and $h_1 \in C$. Hence, $t_1$ has any positive exponent
of $c$ as an initial subword and it follows that $t_1$ is right-attached to
$C$ with respect to $c$ - a contradiction. Thus, we have $ht(C) < ht(h_1)
< ht(w) = ht(g)$ and the claim is proved.

\smallskip

Now, existence of an element $g \in H$ starting with any positive exponent of
$a$ implies existence of such an element $f$ such that $ht(f) =
ht(K)$ which is a contradiction.

\end{proof}

\subsection{Connecting elements}
\label{subs:4.1}

We call a pair of elements $u, v \in CDR(\mathbb{Z}[t],X)$ an
{\em admissible pair} if
\begin{enumerate}
\item $u,v$ are cyclically reduced,
\item $u,v$ are not proper powers,
\item $|u| = |v|$,
\item $u$ is not conjugate to $v^{-1}$ (in particular, $u \neq v^{-1}$).
\end{enumerate}

For an admissible pair $\{u,v\}$ we define an infinite word $s_{u,v} \in
R(\mathbb{Z}[t],X)$, which we call the {\em connecting element} for
the pair $\{u,v\}$, in the following way
\[ \mbox{$s_{u,v}(\beta)$} = \left\{ \begin{array}{ll}
 \mbox{$u(\alpha)$} & \mbox{if $\beta = (k |u| + \alpha, 0), k \geqslant 0,
 1 \leqslant \alpha \leqslant |u|$,} \\
 \mbox{$v(\alpha)$} & \mbox{if $\beta = (-k|v|+\alpha, 1), k \geqslant 1,
 1 \leqslant \alpha \leqslant |v|$.}
\end{array}
\right. \]
Since there exists $m > 0$ such that $u,v \in CDR(\mathbb{Z}^m,X) -
CDR(\mathbb{Z}^{m-1},X)$ then it is easy to see that $s_{u,v} \in
R(\mathbb{Z}^{m+1},X) - R(\mathbb{Z}^m,X)$. Also $s_{u,v}^{-1} =
s_{v^{-1},u^{-1}}$ and $u \circ s_{u,v} = s_{u,v} \circ v$ - both
follow directly from the definition.

Notice that any two connecting elements $s_{u_1,v_1},\ s_{u_2,v_2}$
have the same length whenever $u_1,\ v_1,\ u_2,\ v_2 \in
CDR(\mathbb{Z}^m,X) - CDR(\mathbb{Z}^{m-1},X)$. In this event we have
$$|s_{u_1,v_1}| = |s_{u_2,v_2}| = (0,\ldots,0,1) \in \mathbb{Z}^{m+1}.$$

\begin{lemma}
\label{le:4.1}
Let $u,v$ be elements of a group $H \subset CDR(\mathbb{Z}[t],X)$. If
the pair $\{u,v\}$ is admissible then $s_{u,v} \in CDR(\mathbb{Z}[t],X)$.
\end{lemma}
\begin{proof} Observe that $s_{u,v} \in CDR(\mathbb{Z}[t],X)$ if and only if
$com(s_{u,v}, s_{u,v}^{-1})$ exists.

If $c(u^2, v^{-2}) \geq |u| = |v|$ then $c(s_{u,v}, s_{u,v}^{-1})
\geq |u| = |v|$ then it follows that $u = v^{-1}$ which is impossible
by the assumption. Hence, $c(u^2, v^{-2}) < |u|$, so $c(u^2, v^{-2}) =
c(u,v^{-1})$ and $com(s_{u,v}, s_{u,v}^{-1}) = com(u,v^{-1})$ is defined
since $u,v \in H$.

\end{proof}

\subsection{Main construction}
\label{subs:main_constr}

Now, let $A,B$ be cyclically reduced centralizers in $H$ such that
there exists an isomorphism $\phi : A \rightarrow B$ satisfying the
following conditions
\begin{enumerate}
\item $a$ is not conjugate to $\phi_i(a)^{-1}$ in $H$ for any $a \in A$,
\item $|\phi(a)| = |a|$ for any $a \in A$.
\end{enumerate}

In particular, it follows that $ht(A) = ht(B)$.

\begin{remark}
\label{rem:4.2}
Observe that if $C$ is conjugate to $A$ and $D$ is conjugate to $B$ then
$$\langle H, z \mid z^{-1} A z = B \rangle \simeq \langle H, z' \mid
z'^{-1} C z' = D \rangle.$$
Hence, it is always possible to consider $A$ and $B$ up to taking
conjugates.
\end{remark}

Let $u$ be a generator of $A$ of maximal height, and let $v = \phi(u)
\in B$. Then $v$ is a generator of $B$ of maximal height and $|u| =
|v|$. Observe that from the conditions imposed on $\phi$ it follows that
the pair $u,v$ is admissible. We fix $u$ and $v$ for the rest of the paper.

\smallskip

Now, we are in position to define $s \in R(\mathbb{Z}^{n+1},X)$ which is going to
be an infinite word representing $z$ from the presentation (\ref{eq:G}). By
Lemma \ref{le:attached2} we can assume $A$ to have no right-attached elements
with respect to $u$, and $B$ to have no left-attached elements with respect to
$v$.

Since $\mathbb{Z}[t]$-exponentiation is defined on $CR(\mathbb{Z}[t],X)$
(see \cite{MRS} for details) then for any $f(t) \in \mathbb{Z}[t]$ we can define
$v^{f(t)}, u^{f(t)} \in CR(\mathbb{Z}[t],X)$ so that
$$|v^{f(t)}| = |v| |f(t)|,\ |u^{f(t)}| = |u| |f(t)|$$
and
$$[v^{f(t)},v] = \varepsilon,\ [u^{f(t)},u] = \varepsilon.$$
Thus, if $\alpha = t^{n-ht(A)}$ then $|u^\alpha| = |v^\alpha| = |u| |\alpha|$
and $ht(u^\alpha) = ht(v^\alpha) = ht(u) + (n-ht(A)) = n$. Hence, we
define
$$s = s_{u^\alpha,v^\alpha} \in CDR(\mathbb{Z}^{n+1},X).$$
Observe that $ht(s) = n + 1 = ht(G) + 1$.

\begin{remark}
\label{rem:4.3}
It is easy to see that no element of $H$ has $s^{\pm 1}$ as a subword.
\end{remark}

\begin{lemma}
\label{le:4.2}
For any $h \in A$ we have $s^{-1} \ast h \ast s = \phi(h) \in B$.
\end{lemma}
\begin{proof} By the construction above we have $s^{-1} \ast u \ast s = v = \phi(u)$,
where $u$ is a generator of $A$ of maximal height, and $v$ is a corresponding
generator of $B$ of maximal height. Now, observe that for every $a \in A$
there exists $r > 0$ such that either $a$ or $a^{-1}$ is both an initial and
terminal subword of $u^r$.

So, let us fix $a \in A$ and without loss of generality we assume that
$u^r = a \circ u_1,\ u_1 \in A$. Now,
$$\phi(u_1) \ast \phi(a) = \phi(a) \ast \phi(u_1) = \phi(u^r) = \phi(u)^r = v^r.$$
Since $\phi$ preserves length then
$$|\phi(u_1) \ast \phi(a)| = |\phi(a) \ast \phi(u_1)| = |v^r| = |u^r| = |a| + |u_1| =
|\phi(a)| + |\phi(u_1)|,$$
that is, $\phi(a) \ast \phi(u_1) = \phi(a) \circ \phi(u_1) = \phi(u_1) \circ \phi(a)$.
Hence, the image of an initial subword of $u^r$ of length $\alpha$ representing an
element of $A$ is both an initial and terminal subword of $v^r$ of length $\alpha$.

Finally observe that $s^{-1} \ast a \ast s$ is a terminal subword of $s$ of length
$|a|$.

\end{proof}

Now, our goal is to prove that a pair $H, s$ generates a group in
$CDR(\mathbb{Z}[t],X)$.

\begin{lemma}
\label{le:stab0}
For any $g \in H_i$ there exists $N = N(g) > 0$ such that
$$g \ast u^k = (g \ast u^N) \circ u^{k-N},\ v^k \ast g = v^{k-N} \circ
(v^N \ast g).$$
for any $k > N$.
\end{lemma}
\begin{proof} The required result follows immediately from Lemma
\ref{le:attached3}.

\end{proof}

\begin{lemma}
\label{le:stab1}
\begin{enumerate}
\item[{\bf (i)}] For any $g \in H - A$ there exists $N = N(g) > 0$ such
that for any $k > N$
$$u^{-k} \ast g \ast u^k = u^{-k+N} \circ (u^{-N} \ast g \ast u^N)
\circ u^{k-N}.$$
\item[{\bf (ii)}] For any $g \in H$ there exists $N = N(g) > 0$ such
that for any $k > N$
$$v^k \ast g \ast u^k = v^{k-N} \circ (v^N \ast g \ast u^N)
\circ u^{k-N}.$$
\item[{\bf (iii)}] For any $g \in H - B$ there exists $N = N(g) > 0$ such
that for any $k > N$
$$v^k \ast g \ast v^{-k} = v^{k-N} \circ (v^N \ast g \ast v^{-N})
\circ u^{-k+N}.$$
\end{enumerate}
\end{lemma}
\begin{proof} {\bf (i)} Observe that by Lemma \ref{le:stab0} there exists
$M > 0$ such that $g \ast u^k = (g \ast u^M) \circ u^{k-M}$ for all
$k > M$. So, without loss of generality we can assume $g \ast u = g \circ u$.
Now, if $ht(g) > ht(u)$ then by Lemma \ref{le:attached3} the required result
follows. Now, assume $ht(g) \leq ht(u)$.

Suppose there is no $N > 0$ such that for any $k > N$
$$u^{-k} \ast g \ast u^k = u^{-k+N} \circ (u^{-N} \ast g \ast u^N)
\circ u^{k-N}.$$
Hence for every $K \in \mathbb{N}$ there are $k,m \in \mathbb{N}$ such that
$c(u^k,g \circ u^m) \geqslant M|u|$. It follows that $u^p = g \circ u_2,\ p > 0$,
where $u = u_1 \circ u_2$ and $c((u_1 \circ u_2)^k,(u_2 \circ u_1)^m) \geqslant
2|u|$. Thus, by Lemma \ref{le:LS} it follows that $u_1 \circ u_2 = u_2 \circ
u_1$ and $[g,u] = \varepsilon$ - a contradiction.

\smallskip

{\bf (ii)} Again, by Lemma \ref{le:stab0} there exists $M > 0$ such that
$g \ast u^k = (g \ast u^M) \circ u^{k-M}$ for all $k > M$. So, without loss of
generality we can assume $g \ast u = g \circ u$. Now, if $ht(g) > ht(u)$ then by
Lemma \ref{le:attached3} the required result follows. Now, assume $ht(g) \leq ht(u)$.

Suppose there is no $N > 0$ such that for any $k > N$
$$v^k \ast g \ast u^k = v^{k-N} \circ (v^N \ast g \ast u^N)
\circ u^{k-N}.$$
In other words, for every $K \in \mathbb{N}$ there are $k,m \in \mathbb{N}$ such
that $c(v^{-k},g \circ u^m) \geqslant M|u|$. It follows that $v^{-p} = g \circ v_2,
\ p > 0$, where $v^{-1} = v_1 \circ v_2$ and $u = v_2 \circ v_1 = v_1^{-1} \ast
(v_1 \circ v_2) \ast v_1 = v_1^{-1} \ast v^{-1} \ast v_1$ - a contradiction with
the choice of $\phi$.

\smallskip

{\bf (iii)} The required result follows by the same argument as in {\bf (i)}.

\end{proof}

\begin{defn}
\label{de:s-form}
A sequence
\begin{equation}
\label{eq:p1}
p = (g_1, s^{\epsilon_1}, g_2, \dots, g_k, s^{\epsilon_k}, g_{k+1}),
\end{equation}
where $g_j \in H,\ \epsilon_j \in \{-1,1\}$, $k \geqslant 1$, is called an
{\em $s$-form} over $H$.

An $s$-form (\ref{eq:p1}) is {\em reduced} if subsequences
$$\{s^{-1}, c, s\},\ \ \{s, d, s^{-1}\}$$
where $c \in A,\ d \in B$, do not occur in it.
\end{defn}

Denote by ${\mathcal P}(H,s)$ the set of all $s$-forms over $H$.
We define a partial function $w: {\mathcal P}(H,s) \rightarrow
R(\mathbb{Z}[t],X)$ as follows. If
$$p = (g_1, s^{\epsilon_1}, g_2, \dots, g_k, s^{\epsilon_k}, g_{k+1})$$
then
$$w(p) = (\cdots (g_1 \ast s^{\epsilon_1}) \ast g_2) \ast \cdots \ast g_k)
\ast s^{\epsilon_k}) \ast g_{k+1})$$
if it is defined.

\begin{lemma}
\label{le:p1}
Let $p = (g_1, s^{\epsilon_1}, g_2, \dots, g_k, s^{\epsilon_k}, g_{k+1})$
be an $s$-form over $H$. Then the following hold.
\begin{enumerate}
\item[(1)] The product $w(p)$ is defined and it does not depend on
the placement of parentheses.
\item[(2)] There exists a reduced $s$-form $q$ over $H$ such that $w(q) = w(p)$.
\item[(3)] If $p$ is reduced then there exists a unique representation for $w(p)$ of
the following type
$$w(p) = (g_1 \ast u_1^{N_1}) \circ (u_1^{-N_1} \ast s^{\epsilon_1} \ast
v_1^{-M_1}) \circ (v_1^{M_1} \ast g_2 \ast u_2^{N_2}) \circ \cdots $$
$$\cdots \circ (u_k^{-N_k} \ast s^{\epsilon_k} \ast v_k^{-M_k}) \circ (v_k^{M_k} \ast g_{k+1}),$$
where $N_j, M_j \geq 0,\ u_j = u, v_j = v$ if $\epsilon_j = 1$, and $N_j, M_j \leq 0,\
u_j = v, v_j = u$ if $\epsilon_j = -1$ for $j \in [1,k]$. Moreover, $g_1 \ast u_1^{N_1}$
does not have $u_1^{\pm 1}$ as a terminal subword, $v_{j-1}^{M_{j-1}} \ast g_j \ast
u_j^{N_j}$ does not have $u_j^{\pm 1}$ as a terminal subword for every $j \in [2,k]$,
and $v_{j-1}^{M_{j-1}} \ast g_j \ast s_i^{\epsilon_j}$ does not have $v_{j-1}^{\pm 1}$
as an initial subword for every $j \in [2,k]$, $v_k^{M_k} \ast g_{k+1}$ does not have
$v_k^{\pm 1}$ as an initial subword.
\item[(4)] $w(p) \in CDR(\mathbb{Z}[t],X)$.
\end{enumerate}
\end{lemma}
\begin{proof} Let
$$p = (g_1, s^{\epsilon_1}, g_2, \dots, g_k, s^{\epsilon_k}, g_{k+1})$$
be an $s$-form over $H$.

We show first that (1) implies (2). Suppose that $w(p)$ is defined for every
placement of parentheses and all such products are equal. If $p$ is not reduced
then there exists $j \in [2,k]$ such that either $g_j \in A,\ \epsilon_{j-1}
= -1,\ \epsilon_j = 1$, or $g_j \in B,\ \epsilon_{j-1} = 1,\ \epsilon_j = -1$.
Without loss of generality we can assume the former. Thus, we have
$$s^{-1} \ast g_j \ast s = g'_j \in B \subseteq H$$
and we obtain a new $s$-form
$$p_1 = (g_1, s^{\epsilon_1}, g_2, \dots, g_{j-1} \ast g'_j \ast g_{j+1},
s^{\epsilon_{j+1}}, \dots, g_k, s^{\epsilon_k}, g_{k+1})$$
which is shorter then $p$ and $w(p) = w(p_1)$. Proceeding this way (or by
induction) in a finite number of steps we obtain a reduced $s$-form
$$q = (f_1, s^{\delta_1}, f_2, \ldots, s^{\delta_l}, f_{l+1}),$$
such that $w(q) = w(p)$, as required.

\smallskip

Now we show that (1) implies (3). Assume that
$$p = (g_1, s^{\epsilon_1}, g_2, \dots, g_k, s^{\epsilon_k}, g_{k+1})$$
is reduced.

By Lemma \ref{le:stab0} and Lemma \ref{le:stab1} there exists $r \in
\mathbb{N}$ such that for any $\alpha > r$

\begin{enumerate}
\item[(a)] $g_1 \ast u^\alpha = (g_1 \ast u^r) \circ u^{\alpha-r},\ g_1 \ast
v^{-\alpha} = (g_1 \ast u^{-r}) \circ u^{-\alpha+r}$,
\item[(b)] $v^\alpha \ast g_{k+1} = v^{\alpha-r} \circ (v^r \ast g_{k+1}),\
u^{-\alpha} \ast g_{k+1} = u^{-\alpha+r} \circ (u^{-r} \ast g_{k+1})$,
\item[(c)] $u^{-\alpha} \ast g_j \ast u^\alpha = u^{-(\alpha-r)} \circ (u^{-r} \ast
g_j \ast u^r) \circ u^{\alpha-r}$ for all $j \in [2,k]$,
\item[(d)] $v^\alpha \ast g_j \ast u^\alpha = v^{\alpha-r} \circ (v^r \ast g_j
\ast u^r) \circ u^{\alpha-r}$ for all $j \in [2,k]$,
\item[(e)] $v^\alpha \ast g_j \ast v^{-\alpha} = v^{\alpha-r} \circ (v^r \ast g_j
\ast v^{-r}) \circ u^{-(\alpha-r)}$ for all $j \in [2,k]$.
\end{enumerate}

Hence,
$$w(p) = g_1 \ast s^{\epsilon_1} \ast g_2 \ast \cdots \ast g_k \ast
s^{\epsilon_k} \ast g_{k+1} = $$
$$= (g_1 \ast u_1^r) \circ (u_1^{-r} \ast s^{\epsilon_1} \ast v_1^{-r})
\circ (v_1^r \ast g_2 \ast u_2^r) \circ \cdots \circ (u_k^{-r} \ast
s^{\epsilon_k} \ast v_k^{-r}) \circ (v_k^r \ast g_{k+1}),$$
where $u_j = u, v_j = v$ if $\epsilon_j = 1$ and $u_j = v^{-1}, v_j = u^{-1}$
if $\epsilon_j = -1$ for every $j \in [1,k]$.

Now, if $g_1 \ast u_1^r$ has $u_1^{\gamma_1},\ \gamma_1 \in \mathbb{Z}$ as a
terminal subword then we denote $N_1 = r-\gamma_1$ and rewrite $w(p)$ as follows
$$w(p) = (g_1 \ast u_1^{N_1}) \circ (u_1^{-N_1} \ast s^{\epsilon_1}
\ast v_1^{-r}) \circ (v_1^r \ast g_2 \ast u_2^r) \circ \cdots$$
$$\cdots \circ (u_k^{-r} \ast
s^{\epsilon_k} \ast v_k^{-r}) \circ (v_k^r \ast g_{k+1}).$$
Now, if $v_1^r \ast g_2 \ast u_2^r$ contains $v_1^{\delta_1},\ \delta_1 \in
\mathbb{Z}$ as an initial subword then we denote $M_1 = r-\delta_1$ and again
rewrite $w(p)$
$$w(p) = (g_1 \ast u_1^{N_1}) \circ (u_1^{-N_1} \ast s^{\epsilon_1}
\ast v_1^{-M_1}) \circ (v_1^{M_1} \ast g_2 \ast u_2^r) \circ \cdots$$
$$\cdots \circ (u_k^{-r} \ast s^{\epsilon_k} \ast v_k^{-r}) \circ (v_k^r \ast g_{k+1}).$$
In a finite number of steps we obtain the required result.

\smallskip

Now we prove (1) by induction on $k$. If $k = 1$ then by Lemma \ref{le:stab0}
there exists $r \in \mathbb{N}$ such that
$$g_1 \ast u^\alpha = (g_1 \ast u^r) \circ u^{\alpha-r},\ \ g_1 \ast v^{-\alpha}
= (g_1 \ast v^{-r}) \circ v^{-\alpha+r},$$
$$v^\beta \ast g_2 = u^{\beta-r} \circ (u^r \ast g_2),\ \ u^{-\beta} \ast g_2 =
u^{-\beta+r} \circ (v^{-r} \ast g_2)$$
for any $\alpha, \beta > r$. Hence,
$$(g_1 \ast s^{\epsilon_1}) \ast g_2 = ((g_1 \ast u_1^r) \circ (u_1^{-r} \ast
s^{\epsilon_1})) \ast g_2 =$$
$$= ((g_1 \ast u_1^r) \circ (u_1^{-r} \ast s^{\epsilon_1} \ast v_1^{-r})) \circ
(v_1^r \ast g_2),$$
where $u_1 = u, v_1 = v$ if $\epsilon_1 = 1$ and $u_1 = v^{-1}, v_1 = u^{-1}$ if
$\epsilon_1 = -1$. By Theorem 3.4 \cite{MRS}, the product $(g_1 \ast
s^{\epsilon_1}) \ast  g_2$ does not depend on the placement of
parentheses. So (1) holds for $k = 1$.

Now, consider an initial $s$-subsequence of $p$
$$p_1 = (g_1, s^{\epsilon_1}, g_2, \dots, g_{k-1}, s^{\epsilon_{k-1}}, g_k).$$
By induction $w(p_1)$ is defined and it does not depend on the placement of
parentheses. By the argument above there exists a unique representation of $w(p_1)$
$$w(p_1) = (g_1 \ast u_1^{N_1}) \circ (u_1^{-N_1} \ast s^{\epsilon_1} \ast
v_1^{-M_1}) \circ (v_1^{M_1} \ast g_2 \ast u_2^{N_2}) \circ \cdots$$
$$\cdots \circ (u_{k-1}^{-N_{k-1}} \ast s^{\epsilon_{k-1}} \ast v_{k-1}^{-M_{k-1}}) \circ
(v_{k-1}^{M_{k-1}} \ast g_k),$$
where $N_j, M_j \geq 0,\ u_j = u, v_j = v$ if $\epsilon_j = 1$, and $N_j, M_j \leq 0,\
u_j = v, v_j = u$ if $\epsilon_j = -1$ for $j \in [1,k-1]$. To prove
that $p$ satisfies (1) it suffices to show that
$$w(p_1) \ast (s^{\epsilon_k} \ast g_{k+1})$$
is defined and does not depend on the placement of parentheses.

Without loss of generality we assume $\epsilon_{k-1} = 1, \epsilon_k = 1$ - other
combinations of $\epsilon_{k-1}$ and $\epsilon_k$ are considered similarly.

By Lemma \ref{le:stab1}
$$((u^{-N_{k-1}} \ast s \ast v^{-M_{k-1}}) \circ (v^{M_{k-1}} \ast g_k)) \ast
(s \ast g_{k+1}) = (u^{-N_{k-1}} \ast s \ast v^{-M_{k-1}-r})$$
$$ \circ (v^{M_{k-1}+r} \ast g_k \ast u^{m_1}) \circ (u^{-m_1} \ast s \ast v^{-m_2})
\circ (v^{m_2} \ast g_{k+1})$$
for some $m,r \in \mathbb{N}$. Thus $w(p_1) \ast (s \ast g_{k+1})$
is defined and does not depend on the placement of parentheses.

\smallskip

Now we prove (4). By (3) there exists a unique representation of $w(p)$
$$w(p) = (g_1 \ast u_1^{N_1}) \circ (u_1^{-N_1} \ast s^{\epsilon_1} \ast
v_1^{-M_1}) \circ (v_1^{M_1} \ast g_2 \ast u_2^{N_2}) \circ \cdots$$
$$\cdots \circ (u_k^{-N_k} \ast s^{\epsilon_k} \ast v_k^{-M_k}) \circ (v_k^{M_k} \ast
g_{k+1}),$$
where $N_j, M_j \geq 0,\ u_j = u, v_j = v$ if $\epsilon_j = 1$, and $N_j, M_j
\leq 0,\ u_j = v, v_j = u$ if $\epsilon_j = -1$ for $j \in [1,k]$.
By Lemma 3.8 \cite{MRS}, to prove that $w(q) \in CDR(\mathbb{Z}[t],X)$ it
suffices  to show that
$$g^{-1} \ast w(q) \ast g \in CDR(\mathbb{Z}[t],X)$$
for some $g \in R(\mathbb{Z}[t],X)$.

Without loss of generality we assume $\epsilon_k = 1$ and consider two cases.

\begin{enumerate}

\item[(i)] $g_{k+1} \ast g_1 \ast u_1^{N_1} \notin B$ or $g_{k+1} \ast g_1 \ast u_1^{N_1}
\in B$ but $\epsilon_1 = 1$.

Without loss of generality we assume the former and $\epsilon_1 = 1$. By Lemma
\ref{le:stab1} there exists $N \in \mathbb{N}$ such that
$$(u^{-m} \ast (g_1 \ast u^{N_1})^{-1}) \ast w(p) \ast ((g_1 \ast u^{N_1})
\ast u^m) = (u^{-N_1-m} \ast s \ast v^{-M_1})$$
$$\circ (v^{M_1} \ast g_2 \ast u_2^{N_2}) \circ \cdots \circ (u^{-N_k} \ast s
\ast v^{-M_k-m}) \circ (v^{M_k+m} \ast g_{k+1} \ast g_1 \ast u^{N_1+N}) \circ
u^{m-N}$$
for any $m > N$. Thus,
$$(u^{-m} \ast (g_1 \ast u^{N_1})^{-1}) \ast w(p) \ast ((g_1 \ast u^{N_1})
\ast u^m) \in CR(\mathbb{Z}[t],X) \subset CDR(\mathbb{Z}[t],X).$$

\item[(ii)] $g_{k+1} \ast g_1 \in B$ and $\epsilon_1 = -1$.

Thus, $s \ast (g_{k+1} \ast g_1) \ast s^{-1} \in C_i$ and we have
$$(g_1 \ast s)^{-1} \ast w(p) \ast (g_1 \ast s) =
(g_2 \ast u_2^{N_2}) \circ (u_2^{-N_2} \ast s^{\epsilon_2} \ast
v_2^{-M_2}) \circ \cdots$$
$$\cdots \circ (u_{k-1}^{-N_{k-1}} \ast s^{\epsilon_{k-1}}
\ast v_{k-1}^{-M_{k-1}}) \circ (v_{k-1}^{M_{k-1}} \ast g_k),$$
so the number of $s^{\pm 1}$ is reduced by two and we can use induction.
\end{enumerate}

\end{proof}

Now we are ready to prove the main results of this subsection from which
Theorem \ref{main4} follows.

\begin{theorem}
\label{th:U(P)}
Put
$$P = P(H,s) = \{ g \ast s^\epsilon \ast h \mid g,h \in H, \epsilon
\in \{-1,0,1\} \} \subseteq CDR(\mathbb{Z}[t],X).$$
Then the following hold.
\begin{itemize}
\item [(1)] $P$ generates a subgroup $H^*$ in $CDR(\mathbb{Z}[t],X)$.
\item [(2)] $P$, with the multiplication $\ast$ induced from $R(\mathbb{Z}[t],X)$, is a pregroup and
$H^*$ is isomorphic to $U(P)$.
\item [(3)] $H^*$ is isomorphic to $G = \langle H, z \mid z^{-1} A z = B \rangle$.
\end{itemize}
\end{theorem}
\begin{proof} We need the following claims.

\medskip

{\bf Claim 1.} Let $g_j \ast s^{\epsilon_j} \ast h_j  \in P$ $j = 1,2$. If
$$g_1 \ast s^{\epsilon_1} \ast h_1 = g_2 \ast s^{\epsilon_2} \ast
h_2$$
then $\epsilon_1 = \epsilon_2$ and $h_1 \ast h_2^{-1} \in A$ if
$\epsilon_1 = -1$, and $h_1 \ast h_2^{-1} \in B$ if
$\epsilon_1 = 1$.

\smallskip

To prove the claim consider an $s$-form
$$a = (g_1, s^{\epsilon_1}, h_1 \ast h_2^{-1}, s^{-\epsilon_2},
g_2^{-1}).$$
By Lemma \ref{le:p1}, $w(a)$ is defined and
$$g_1 \ast s^{\epsilon_1} \ast h_1 \ast h_2^{-1} \ast s^{-\epsilon_2}
\ast g_2^{-1} = \varepsilon.$$
Hence, $a$ is not reduced and the claim follows.

\medskip

For every $p \in P$ we fix now a representation $p = g_p \ast
s^{\epsilon_p} \ast h_p$, where $g_p, h_p \in H, \epsilon_p \in
\{-1,0,1\}$.

\medskip

{\bf Claim 2.} Let $p = g_p \ast s^{\epsilon_p} \ast h_p, \ \ q =
g_q \ast s^{\epsilon_q} \ast h_q$ be in $P$. If $p \ast q \in P$
then either $\epsilon_p \epsilon_q = 0$, or
$\epsilon_p = - \epsilon_q \neq 0$ and $h_p \ast g_q \in A$ if
$\epsilon_p = -1$, and $h_p \ast g_q \in B$ if $\epsilon_q = 1$.

\smallskip

Let $x = p \ast q \in P$ and $x = g_x \ast s^{\epsilon_x} \ast h_x$.
Assume that $\epsilon_p \epsilon_q \neq 0$.

\smallskip

\begin{enumerate}

\item[(a)] $\epsilon_x \neq 0$

\smallskip

Consider an $s$-form
$$a = (g_p, s^{\epsilon_p}, h_p \ast g_q, s^{\epsilon_q}, h_q \ast
g_x, s^{\epsilon_x}, h_x).$$
By Lemma \ref{le:p1}, $w(a)$ is defined and
$$w(a) = g_p \ast s^{\epsilon_p} \ast h_p \ast g_q \ast s^{\epsilon_q}
\ast h_q \ast g_x \ast s^{\epsilon_x} \ast h_x = \varepsilon.$$
Hence, $a$ is not reduced and either a subsequence
$$\{s^{\epsilon_p}, h_p \ast g_q, s^{\epsilon_q}\},$$
or a subsequence
$$\{s^{\epsilon_q}, h_q \ast g_x, s^{\epsilon_x}\}$$
is reducible. In the former case we are done, so assume that
$\{s^{\epsilon_q}, h_q \ast g_x, s^{\epsilon_x}\}$ can be reduced.
Without loss of generality we can assume that $\epsilon_q = -1,\
\epsilon_x = 1,\ h_q \ast g_x \in A$. Hence,
$$s^{\epsilon_q} \ast h_q \ast g_x \ast s^{\epsilon_x} = g \in B$$
and we have
$$w(a) = g_p \ast s^{\epsilon_p} \ast h_p \ast g_q \ast g \ast h_x =
\varepsilon.$$
Now, it follows $\epsilon_p = 0$ - a contradiction with our assumption.

\smallskip

\item[(b)] $\epsilon_x = 0$

\smallskip

Hence, $x = g \in H$ and we consider an $s$-form
$$a = (g_p, s^{\epsilon_p}, h_p \ast g_q, s^{\epsilon_q}, h_q \ast g).$$
By Lemma \ref{le:p1}, $w(a)$ is defined and
$$w(a) = g_p \ast s^{\epsilon_p} \ast h_p \ast g_q \ast s^{\epsilon_q}
\ast h_q \ast g = \varepsilon.$$
\end{enumerate}

Now, the claim follows automatically.

\medskip

Below we call a tuple $y = (y_1, \ldots, y_k) \in P^k$ a {\em reduced} $P$-sequence if
$y_j \ast y_{j+1} \notin P$ for $j \in [1,k-1]$. Observe, that if $y = (y_1, \ldots, y_k)$ is a
reduced $P$-sequence and $y_j = g_j \ast s_i^{\epsilon_j} \ast h_j$ then either $k \leq 1$ or
$y$ has the following properties which follow from Claim 2:
\begin{enumerate}
\item[(a)] $\epsilon_j \neq 0$ for all $j \in [1,k]$,
\item[(b)] if $\epsilon_j = -1,\ \epsilon_{j+1} = 1$ then $h_j \ast g_{j+1} \notin A$ for $j \in [1,k-1]$,
\item[(c)] if $\epsilon_j = 1,\ \epsilon_{j+1} = -1$ then $h_j \ast g_{j+1} \notin B$ for $j \in [1,k-1]$.
\end{enumerate}
In particular, the $s$-form over $H$
$$p_y = (g_1,\ s^{\epsilon_1},\ h_1 \ast g_2,\ s^{\epsilon_2},\
\ldots,\ h_{n-1} \ast g_n,\ s^{\epsilon_k},\ h_n),$$
is reduced.

\medskip

To prove (1) observe first that $P^{-1} = P$. Now if $y_1, \ldots, y_k \in
P$ then $y_1 \ast \cdots \ast y_k = w(p_y)$, where $y = (y_1, \ldots, y_k)$. Hence, by Lemma \ref{le:p1},
the product $y_1 \ast y_2 \ast \cdots \ast y_k$ is defined in $CDR(\mathbb{Z}[t],X)$ and it belongs to
$CDR(\mathbb{Z}[t],X)$. It follows that $H^* = \langle P \rangle$ is a subgroup of $CDR(\mathbb{Z}[t],X)$
which consists of all words $w(p)$, where $p$ ranges through all possible $s$-forms over $H$. Hence, (1)
is proved.

\smallskip

Now we prove (2). By Theorem 2, \cite{R2}, to prove that $P$ is a pregroup and the inclusion $P \to H^*$ extends to an
isomorphism $U(P) \simeq H^*$ it is enough to show that all reduced $P$-sequences representing the same element
have the same $P$-length.
\smallskip

Suppose two reduced $P$-sequences
$$(u_1, u_2, \ldots, u_k),\ (v_1, v_2, \ldots, v_n)$$
represent the same element $g \in H^*$. That is,
$$(u_1 \ast \cdots \ast u_k) \ast (v_1 \ast \cdots \ast v_n)^{-1} = \varepsilon.$$
We use the induction on $k+n$ to show that $k = n$. Observe that $k = 0$ implies
$n = 0$, otherwise we get a contradiction with Lemma \ref{le:p1} (3). Hence, we can assume
$k, n > 0$, that is, $k + n \geq 2$. If the $P$-sequence
$$a = (u_1, \ldots, u_k, v_n^{-1}, \ldots, v_1^{-1})$$
is reduced then the underlying $s$-form is reduced and hence, by Lemma \ref{le:p1} (3)
$$w(a) = u_1 \ast \ldots \ast u_k \ast v_n^{-1} \ast \ldots \ast v_1^{-1} \neq \varepsilon.$$
Hence,
$$(u_1, \ldots, u_k, v_n^{-1}, \ldots, v_1^{-1})$$
is not reduced and $u_k \ast v_n^{-1} \in P$. If $u_k = g_1 \ast s^{\epsilon_1} \ast h_1,\
v_n = g_2 \ast s^{\epsilon_2} \ast h_2$, where $g_i, h_i \in H$ and $\epsilon_i \in
\{-1,0,1\},\ i = 1,2$ then by  Claim 2 either $\epsilon_1 \epsilon_2 = 0$, or
$\epsilon_1 = \epsilon_2 \neq 0$ and $h_1 \ast h_2^{-1} \in A$ if
$\epsilon_1 = -1$, and $h_1 \ast h_2^{-1} \in B$ if $\epsilon_1 = 1$. In the former case, for
example, if $\epsilon_2 = 0$ then $n = 1,\ v_n \in H$
and $b = (u_1, \ldots, u_k \ast v_n^{-1})$ is a reduced $P$-sequence such that $w(b) = \varepsilon$
- a contradiction with Lemma \ref{le:p1} (3) unless $k = 1,\ u_1 \in H$. In the latter case,
$u_k \ast v_n^{-1} \in H$ and it follows that
$$(u_1, u_2, \ldots, u_{k-1} \ast (u_k \ast v_n^{-1})),\ (v_1, v_2, \ldots, v_{n-1})$$
represent the same element in $H^*$ while the sum of their lengths is less
than $k+n$, so the result follows by induction.

\medskip

Finally, to prove (3) observe first that $H$ embeds into $G$. We denote this embedding
by $\theta$. Now we define a map $\phi : P \rightarrow G$ as
follows. For $g \ast s^\epsilon \ast h \in P$ put
$$g \ast s^\epsilon \ast h \ \stackrel{\phi}{\rightarrow}\
\theta(g)\ z^\epsilon\ \theta(h).$$
It follows from Claim 2 that $\phi$ is a morphism of pregroups.
Since $H^* \simeq U(P)$, the morphism
$\phi$ extends to a unique homomorphism $\psi : H^* \rightarrow
G$. We claim that $\psi$ is bijective. Indeed, observe first that $G =
\langle H, z \rangle$. Now, since $\psi(s^\epsilon) = z^\epsilon$ and
$\psi = \phi = \theta$ on $H$, it follows that $\psi$ is onto. To see
that $\psi$ is one-to-one it suffices to notice that if
$$(g_1 \ast s^{\epsilon_1} \ast h_1, g_2 \ast s^{\epsilon_2}
\ast h_2, \ldots, g_m \ast s^{\epsilon_m} \ast h_m)$$
is a reduced $P$-sequence then
$$y = (g_1, s^{\epsilon_1}, h_1 \ast g_2, s^{\epsilon_2}, \ldots, s^{\epsilon_m}, h_m)$$
is a reduced $s$-form and $w(y)^\psi \neq 1$ by Britton's
Lemma (see, for example, \cite{MKS0}). This proves that $\psi$ is an
isomorphism, as required.

\end{proof}

\begin{theorem}
\label{th:Greg} Let $G = \langle H, z \mid z^{-1} A z = B \rangle$. Then, in the notation above,  the free length function on $L : G \rightarrow \mathbb{Z}^{n+1}$ induced by the
isomorphism $\psi : H^* \rightarrow G$ is regular.
\end{theorem}
\begin{proof} Observe that is is enough to show that the length function induced
on $H^* = \langle P \rangle$ from $CDR(\mathbb{Z}[t],X)$ is regular.

Let $g,h \in H^*$. Then $g$ and $h$ can be written in the unique normal forms
$$g = g_1 \circ (u_1^{-N_1} \ast s^{\epsilon_1} \ast v_1^{-M_1}) \circ g_2 \circ \cdots \circ
(u_k^{-N_k} \ast s^{\epsilon_k} \ast v_k^{-M_k}) \circ g_{k+1},$$
$$h = h_1 \circ (w_1^{-L_1} \ast s^{\delta_1} \ast x_1^{-P_1}) \circ h_2 \circ \cdots \circ
(w_m^{-L_m} \ast s^{\delta_m} \ast x_m^{-P_m}) \circ h_{m+1},$$
where $N_j, M_j \geq 0,\ u_j = u, v_j = v$ if $\epsilon_j = 1$, and $N_j, M_j \leq 0,\
u_j = v, v_j = u$ if $\epsilon_j = -1$ for $j \in [1,k];\ L_i, P_i \geq 0,\ w_i = u, x_i = v$
if $\delta_i = 1$, and $L_i, P_i \leq 0,\ w_i = v, x_i = u$ if $\delta_i = -1$ for $i \in [1,m]$.
Moreover, $g_1$ does not have $u_1^{\pm 1}$ as a terminal subword, $g_j$ does not have
$u_j^{\pm 1}$ as a terminal subword for every $j \in [2,k]$, and $g_j$ does not have $v_{j-1}^{\pm 1}$
as an initial subword for every $j \in [2,k]$, $g_{k+1}$ does not have
$v_k^{\pm 1}$ as an initial subword; $h_1$ does not have $w_1^{\pm 1}$ as a terminal subword,
$h_i$ does not have $w_i^{\pm 1}$ as a terminal subword for every $i \in [2,m]$, and $h_i$ does not have
$x_{j-1}^{\pm 1}$ as an initial subword for every $i \in [2,m]$, $h_{m+1}$ does not have
$x_m^{\pm 1}$ as an initial subword.

\smallskip

If there exist $k_1,k_2 > 0$ such that
$$c = c(g,h) \leq \min\{|g_1 \circ u_1^{k_1}|,|h_1 \circ w_1^{k_2}|\}$$
then $com(g,h) = com(g_1 \circ u_1^{k_1}, h_1 \circ w_1^{k_2}) \in H$. Now, assume
that $r \in [1,k]$ is the minimal natural number such that $c$ is an initial subword of
$$f_1 = g_1 \circ (u_1^{-N_1} \ast s^{\epsilon_1} \ast v_1^{-M_1}) \circ g_2 \circ \cdots \circ
(u_r^{-N_r} \ast s^{\epsilon_r} \ast v_r^{-M_r}) \circ g_{r+1} \circ u_r^{p_1},$$
where $p_1 \in \mathbb{Z}$. Similarly, assume that $q \in [1,m]$ is the minimal natural
number such that $c$ is an initial subword of
$$f_2 = h_1 \circ (w_1^{-L_1} \ast s^{\delta_1} \ast x_1^{-P_1}) \circ h_2 \circ \cdots \circ
(w_q^{-L_q} \ast s^{\delta_q} \ast x_q^{-P_q}) \circ h_{q+1} \circ w_q^{p_2},$$
where $p_2 \in \mathbb{Z}$.
From uniqueness of normal forms it follows that $r = q$ and we have $g_i = h_i,\ u_i =
w_i,\ v_i = x_i,\ N_i = L_i,\ \epsilon_i = \delta_i,\ i \in [1,r]$ and $M_i = P_i,\ i
\in [1,r-1]$.

\smallskip

Without loss of generality we can assume $\epsilon_r = 1$. Hence, $v_r = x_r = v$.

Observe that $c$ can be represented as a concatenation $c = c_1 \circ c_2$, where
$$c_1 = g_1 \circ (u_1^{-N_1} \ast s^{\epsilon_1} \ast v_1^{-M_1}) \circ g_2 \circ \cdots \circ
(u^{-N_r} \ast s \ast v^{-l})$$
and $l \geq \max\{M_r,P_r\}$, and
$$c_2 = com(v^{l-M_r} \circ g_{r+1} \circ u_r^{p_1}, v^{l-P_r} \circ h_{r+1} \circ w_r^{p_2}).$$
Obviously, $c_1 \in H^*$. Also, $c_2 \in H$ since $v^{l-M_r} \circ g_{r+1} \circ u_r^{p_1},\
v^{l-P_r} \circ h_{r+1} \circ w_r^{p_2} \in H$ and the length function on $H$ is regular.
Hence, $c \in H^*$.

\end{proof}


\end{document}